\documentclass[a4paper,reqno]{amsart}

\usepackage[utf8]{inputenc}
\usepackage[english]{babel}

\usepackage{hyperref} 
\hypersetup{hidelinks,breaklinks}

\usepackage{amsmath}
\usepackage{amssymb}
\usepackage{amsthm} 
\usepackage{thmtools} 
\usepackage{mathtools} 

\usepackage{tikz} 
\usetikzlibrary{cd}

\usepackage{multirow}

\usepackage{dsfont} 
\usepackage{pifont} 
\usepackage{nicefrac} 

\usepackage{ifthen}
\usepackage{enumitem}

\newcounter{i} 
\newtoks\striche 

\newcounter{Resultate}[section] 
\renewcommand{\theResultate}{\thesection.\arabic{Resultate}} 

\numberwithin{equation}{section} 

\newcommand{\myref}[2]{\hyperref[#1]{#2 \ref*{#1}}}

\newcommand{\R}{\mathbb{R}}
\newcommand{\N}{\mathbb{N}} 

\newcommand{\C}{\mathbb{C}}
\newcommand{\K}{\mathbb{K}}

\newcommand{\mc}{\mathcal}

\newcommand{\la}{\lambda}

\newcommand{\Lp}[1]{L^{#1}} 
\newcommand{\Lb}[2][]{\mathcal{L}_{\mathrm{b}}#1(#2#1)} 

\newcommand{\indicator}{\mathds{1}} 
\newcommand{\diffd}{\mathrm{d}} 
\newcommand{\dx}[1][x]{\,\diffd#1}

\renewcommand{\vec}[1]{\mathbf{#1}}
\newcommand{\myvec}[1]{\boldsymbol{#1}} 
\newcommand{\cl}[2][]{\overline{#2}\ifthenelse{ \equal{#1}{} }{}{^{#1}}} 
\newcommand{\quspace}[3][]{{\raisebox{.2em}{$#1#2$}\mspace{-4.5mu}\left/\mspace{2mu}\raisebox{-.2em}{$#1#3$}\right.}} 

\DeclareMathOperator{\ran}{ran}
\DeclareMathOperator{\mul}{mul}

\DeclareMathOperator{\supp}{supp}
\DeclareMathOperator{\dom}{dom}

\DeclareMathOperator{\Div}{div}
\DeclareMathOperator{\grad}{grad}
\DeclareMathOperator{\rot}{rot}

\DeclarePairedDelimiter{\set}{\{}{\}}
\DeclarePairedDelimiter{\norm}{\lVert}{\rVert}
\DeclarePairedDelimiter{\abs}{\vert}{\vert}
\DeclarePairedDelimiter{\scprod}{\langle}{\rangle}
\DeclarePairedDelimiter{\dualprod}{\langle}{\rangle}

\renewcommand{\Re}{\operatorname{Re}}

\newcommand{\dual}{'}

\newcommand{\adjunsymb}{\ast} 
\newcommand{\adjunX}[1]{^{\adjunsymb_{#1}}}

\newcommand{\adjun}[1][1]{%
  \setcounter{i}{1}%
  \striche={\adjunsymb}%
  \loop%
  \ifnum\value{i}<#1%
  \striche=\expandafter{\the\expandafter\striche\adjunsymb}%
  \stepcounter{i}%
  \repeat%
  ^{\the\striche}%
}

\newcommand{\mapping}[4]{%
  \left\{%
    \begin{array}{rcl}%
      #1	&\to		& #2,		\\
      #3 	&\mapsto 	& #4
    \end{array}%
  \right.%
}

\newcommand{\Zint}[2]{\set{#1,\dots,#2}}
\newcommand{\hermitian}{^{\mathsf{H}}}
\newcommand{\trans}{^{\mathsf{T}}}

\newcommand{\hamiltonian}{\mathcal{H}}
\newcommand{\idop}{\mathrm{I}}
\newcommand{\opid}{\idop}
\newcommand{\myL}{L}
\newcommand{\diffop}[1][]{\ifthenelse{\equal{#1}{}}{\myL}{#1}_{\mspace{-2.5mu}\partial}}
\newcommand{\diffopad}{\myL\hermitian_{\mspace{-2.5mu}\partial}}
\newcommand{\Lnu}[1][]{\ifthenelse{\equal{#1}{}}{}{\indicator_{#1}\mspace{-2.5mu}}\bar{\myL}_{\mspace{-2.5mu}\nu}}
\newcommand{\Lnuad}{\bar{\myL}\hermitian_{\mspace{-2.5mu}\nu}}
\newcommand{\Lnub}{L_{\mspace{-2.5mu}\nu}}
\newcommand{\Lnubad}{L\hermitian_{\mspace{-2.5mu}\nu}}

\newcommand{\piL}[1][]{\pi_{\myL} \ifthenelse{\equal{#1}{}}{}{^{#1}} }

\newcommand{\myP}{P}
\newcommand{\Pop}{\diffop[\myP]}
\newcommand{\Pnu}{\bar{\myP}_{\mspace{-2.5mu}\nu}}

\newcommand{\Pspace}{H\mspace{-1.5mu}(\Pop,\Omega)}

\newcommand{\Lspace}[1][]{H \ifthenelse{\equal{#1}{}}{}{_{\mspace{-0.7mu}#1}\mspace{1.5mu}}\mspace{-1.5mu}(\diffop,\Omega)}
\newcommand{\Ladspace}[1][]{H \ifthenelse{\equal{#1}{}}{}{_{\mspace{-0.7mu}#1}\mspace{1.5mu}}\mspace{-1.5mu}(\diffopad,\Omega)}
\newcommand{\Lnorm}[2][]{\norm[#1]{#2}_{\mspace{-1.5mu}\Lspace}}
\newcommand{\Ladnorm}[2][]{\norm[#1]{#2}_{\mspace{-1.5mu}\Ladspace}}
\newcommand{\PLspace}[1][]{\ifthenelse{\equal{#1}{}}{\mathcal{V}_{\myL}}{\mathcal{V}_{\myL,#1}}} 
\newcommand{\PLadspace}{\mathcal{V}_{\myL\hermitian}} 

\newcommand{\nrg}{E} 
\newcommand{\XH}{\mathcal{X}_{\hamiltonian}}

\newcommand{\boundtr}{\gamma_{0}}

\newcommand{\irow}{m_{1}}
\newcommand{\icol}{m_{2}}

\newcommand{\hs}{\hssymbol_{+}}
\newcommand{\hsmid}{\hssymbol_{0}}
\newcommand{\hsdual}{\hssymbol_{-}}

\theoremstyle{plain}
\newtheorem{theorem}{Theorem}[section]
\newtheorem{lemma}[theorem]{Lemma}
\newtheorem{proposition}[theorem]{Pro\-po\-si\-tion}
\newtheorem{corollary}[theorem]{Corollary}

\theoremstyle{definition}
\newtheorem{definition}[theorem]{Def\-ini\-tion}
\newtheorem{example}[theorem]{Example}
\newtheorem{assumption}[theorem]{Assumption}

\theoremstyle{remark}
\newtheorem{remark}[theorem]{Remark}

\newenvironment{metaUmgebung}[4][]{%
	\par\addvspace{\medskipamount}%
	\noindent%
	\refstepcounter{Resultate}%
		{#3#2 \theResultate{}.}%
		\ifthenelse{\equal{#1}{}}{}{ (#1)}%
		#4\space\ignorespaces%
	}%
	{
	\par\addvspace{\medskipamount}%
	}%

\newenvironment{Definition}{%
	\begin{definition}%
	}%
	{\end{definition}%
	}
	
\newenvironment{Theorem}[1][]{%
	\begin{theorem}%
	}%
	{\end{theorem}
	}

\newenvironment{Corollary}[1][]{%
	\begin{corollary}
	}%
	{%
	\end{corollary}
	}

\newenvironment{Proposition}[1][]{%
	\begin{proposition}
	}%
	{%
	\end{proposition}
	}
	
\newenvironment{Lemma}[1][]{%
	\begin{lemma}%
	}%
	{%
	\end{lemma}
	}

	{%
	\end{metaUmgebung}%
	}

\newenvironment{Remark}{%
	\begin{remark}%
	}%
	{
	\end{remark}%
	}
	
\newenvironment{Example}{%
	\begin{example}
	}%
	{%
	\end{example}
	}

\newenvironment{Proof}[1][]{\ifthenelse{\equal{#1}{}}{\begin{proof}}{\begin{proof}[#1]}}{\end{proof}}

\newenvironment{Voraussetzung}{%
	\begin{assumption}%
	}%
	{%
	\end{assumption}%
	}

\title[Port-Hamiltonian Systems on multidimensional spatial domains]{Well-posedness of linear first order Port-Hamiltonian Systems on multidimensional spatial domains}
\author{Nathanael Skrepek}
\thanks{This project has received funding from the European Union’s Horizon 2020 research and innovation programme under the Marie Sklodowska-Curie grant agreement No 765579}
\date{\today}
\address{University of Wuppertal \\
  School of Mathematics and Natural Science \\
  Gau{\ss}stra{\ss}e 20 \\
  D-42119 Wuppertal \\
  Germany}
\email{skrepek@uni-wuppertal.de}
\begin{document}

\maketitle

\begin{abstract}
  We consider a port-Hamiltonian system on a spatial domain $\Omega \subseteq \R^n$ that is bounded with Lipschitz boundary.
  We show that there is a boundary triple associated to this system. Hence, we can characterize all boundary conditions that provide unique solutions that are non-increasing in the Hamiltonian.
  As a by-product we develop the theory of quasi Gelfand triples.
  Adding ``natural'' boundary controls and boundary observations yields scattering/impedance passive boundary control systems.
  This framework can be applied to the wave equation, Maxwell equations and Mindlin plate model, and probably many more.
\end{abstract}

\section{Introduction}

The aim of this paper is to develop a port-Hamiltonian framework on multidimensional spatial domains that justifies existence and uniqueness of solutions.
Those systems can be described by the following equations
\begin{align*}
\frac{\partial}{\partial t} x(t,\zeta)	
  &= \sum_{i=1}^{n} \frac{\partial}{\partial \zeta_{i}} P_{i} \big(\hamiltonian(\zeta) x(t, \zeta)\big)
+ P_{0}\big(\hamiltonian(\zeta)x(t,\zeta)\big),
  && \zeta \in \Omega, t \geq 0, \\
  x(0,\zeta) &= x_{0}(\zeta),
  && \zeta \in \Omega,
\end{align*}
where $P_{i}$ and $P_{0}$ are matrices, $\hamiltonian$ is the Hamiltonian density, and $\Omega$ is a open subset of $\R^{n}$ with bounded Lipschitz boundary.
We will restrict ourselves to the case, where the matrices $P_i$ have the block shape $\left[\begin{smallmatrix}0 & L_i \\ L_i\hermitian & 0\end{smallmatrix}\right]$ for $i\in\Zint{1}{n}$.
We also introduce ``natural'' boundary controls and observations which makes the system a scattering passive (engery preserving) or impedance passive (energy preserving) boundary control system.

The port-Hamiltonian formulation has proven to be a powerful tool for the modeling and control of complex mutliphysics systems. An introductory overview can be found in \cite{port-hamiltonian-overview}. For a one-dimensional spatial domain concerns about existence and uniqueness of solutions are covered in \cite{jacob-zwart-book}.

Chapter 8 of the Ph.D. thesis \cite{phd-villegas} also regards such port-Hamiltonian systems that have multidimensional spatial domains, but the results demand very strong assumptions, which are in case of the Maxwell equations and Mindlin plate model not satisfied.
With the following approach we will overcome these limits.

The strategy is to find a boundary triple associated to the differential operator.
The multidimensional integration by parts formula already suggests possible operators for a boundary triple, but unfortunately these operators cannot be extended to the entire domain of the differential operator.
Hence, we need to adapt the codomain of these boundary operators, which will lead to the construction of suitable boundary spaces for this problem.
These boundary spaces behave like a Gelfand triple with the original codomain as pivot space, but lack of a chain inclusion.

Up to the authors best knowledge there is no theory about this setting.
So we will develop the notion of \hyperref[sec:quasi-gelfand-triple]{\emph{quasi Gelfand triples}} in \autoref{sec:quasi-gelfand-triple}, which equips us with the tools to state the boundary condition in terms of the pivot space instead of the artificially constructed boundary spaces (\autoref{th:key-theorem}).

The approach to the wave equation in \cite{kurula-zwart-wave} perfectly fits the framework presented in this paper.
In fact, many ideas from \cite{kurula-zwart-wave} are generalized in this work. 
Also the Maxwell equations can be formulated as such a port-Hamiltonian system and the results in \cite{weiss-staffans-maxwell} can also be derived with the tools of this paper.
Moreover, this theory can be applied on the model of Mindlin Plate in \cite{mindlin-plate-andrea},\cite{mindlin-plate}.
In \autoref{sec:boundary-control-system} we give examples of how this framework can be applied to these three PDEs.

\section{Boundary Triple}

In this section we state the most important properties of boundary triples for skew-symmetric operators for this work.
More details can be found in \cite[chapter 3.4]{gorbachuk} and \cite{kurula-zwart-wave}.

A linear relation $T$ between two vector spaces $X$ and $Y$ is a linear subspace of $X\times Y$. Clearly, every linear operator is also a linear relation. We will use the following notation
\begin{align*}
  \ker T &:= \set{x \in X: (x,0) \in T},  &\ran T &:= \set{y \in Y: \exists x :(x,y)\in T}, \\
  \mul T &:= \set{y \in Y: (0,y) \in T},  &\dom T &:= \set{x \in X: \exists y :(x,y) \in T}.
\end{align*}
Thus, $T$ is single-valued, if $\mul T = \set{0}$.
For a linear relation $T$ between two Hilbert spaces $X$ and $Y$ the adjoint relation is defined by
\[
  T\adjun := \set{(u,v)\in Y \times X : \scprod{u,y}_{Y} = \scprod{v,x}_{X} \;\text{for all}\; (x,y) \in T}
\]
and the following holds true
\[
  \ker T\adjun = (\ran T)^{\perp}
  ,\quad
  \mul T\adjun = (\dom T)^{\perp}
  \quad\text{and}\quad
    T\adjun = \left[\begin{smallmatrix} 0 & \opid_{Y} \\ -\opid_{X} & 0\end{smallmatrix}\right] T^{\perp}
  ,
\]
where $\left[\begin{smallmatrix} 0 & \opid_{Y} \\ -\opid_{X} & 0\end{smallmatrix}\right] T := \set{(y,-x): (x,y)\in T}$.
A linear relation $T$ on a Hilbert space $H$ (between $H$ and $H$) is \emph{dissipative}, if $\Re \scprod{x,y}_{H} \leq 0$ for every $(x,y)\in T$ and \emph{maximal dissipative}, if additionally there is no proper dissipative extension of $T$. More details can be found in \cite{linear-relations-book}.

\begin{Definition}\label{def:boundary-triple}
  Let $A_{0}$ be a densely defined, skew-symmetric, and closed operator on a Hilbert space $X$.
  By a \emph{boundary triple} for $A_{0}\adjun$ we mean a triple $(\mc B, B_1, B_2)$ consisting of a Hilbert space $\mc B$, and two linear operators $B_1,B_2 : \dom A_{0}\adjun \to \mc B$ such that
  \begin{enumerate}[label = \textrm{\textup{(\roman*)}}]
  \item\label{def:boundary-triple-surjective} the mapping
  $\big[\begin{smallmatrix}B_1 \\ B_2 \end{smallmatrix}\big]: \dom A_{0}\adjun \to \mc B \times \mc B, x \mapsto \big[\begin{smallmatrix}B_1 x \\ B_2 x\end{smallmatrix}\big]$
    is surjective, and

  \item\label{def:boundary-triple-equation} for $x,y \in \dom A_{0}\adjun$ there holds
    \begin{equation}\label{eq:boundary-triple}
      \scprod{A_{0}\adjun x, y}_{X} + \scprod{x, A_{0}\adjun y}_{X} = \scprod{B_1 x, B_2 y}_{\mc B} + \scprod{B_2 x, B_1 y}_{\mc B}.
    \end{equation}
    
  \end{enumerate}
\end{Definition}

The operator $A_{0}$ can be restored from by restricting $-A_{0}\adjun$ to $\ker B_{1} \cap \ker B_{2}$ as the next lemma will show. However, if $A_0\adjun$ wasn't the adjoint of a skew-symmetric operator then this would not hold as \autoref{ex:b-triple-counter-ex} demonstrates.

\begin{Lemma}\label{le:reconstruct-A0}
  Let $A_{0}$ be a densely defined, skew-symmetric, and closed operator on a Hilbert space $X$ and $(\mc B, B_1, B_2)$ be a boundary triple for $A_0\adjun$.
  Then $A_0 = -A_0\adjun\big\vert_{\ker B_1 \cap \ker B_2}$.
\end{Lemma}

\begin{Proof}
  Let $x \in \ker B_1 \cap \ker B_2$ and $y \in \dom A_0\adjun$. Then the right-hand-side of \eqref{eq:boundary-triple} is $0$. Hence,
  \begin{equation*}
    \scprod{x,A_0\adjun y}_X = \scprod{-A_0\adjun x,y}_X \quad \text{for all}\quad y\in \dom A_0\adjun
    .
  \end{equation*}
  This yields $(x,-A_0\adjun x) \in A_0\adjun[2] = A_0$. Hence, $-A_0\adjun \big\vert_{\ker B_1\cap \ker B_2} \subseteq A_0$.

  On the other hand if $x\in \dom A_0$, then $A_0\adjun x = -A_0 x$ and consequentely
  \begin{equation*}
    \scprod{A_0\adjun x,y}_{X} + \scprod{x,A_0\adjun y}_{X} = \scprod{-x,A_0\adjun y}_{X} + \scprod{x,A_0\adjun y}_{X} = 0
    .
  \end{equation*}
  Therefore, using \eqref{eq:boundary-triple} yields $\scprod*{ \left[\begin{smallmatrix}B_1 \\ B_2\end{smallmatrix}\right] x,\left[\begin{smallmatrix}B_2 \\ B_1\end{smallmatrix}\right] y}_{\mc B \times \mc B} = 0$ for all $y \in \dom A_{0}\adjun$. Since $\left[\begin{smallmatrix}B_2 \\ B_1\end{smallmatrix}\right]$ is surjective on $\mc B \times \mc B$, we have
  \begin{equation*}
    \begin{bmatrix}
      B_1 x \\ B_2 x
    \end{bmatrix}
    \perp
    \mc B \times \mc B
    ,
  \end{equation*}
  which yields $x \in \ker B_1 \cap \ker B_2$.
\end{Proof}

The following result is Theorem 2.2 from \cite{kurula-zwart-wave}.

\begin{Proposition}\label{pr:boundary-triple-facts}%
  \newcommand{\bvec}[1]{\begin{bmatrix}#1\end{bmatrix}}%
  Let $A_{0}$ be a skew-symmetric operator and $(\mathcal{B},B_{1},B_{2})$ be a boundary triple for $A_{0}\adjun$. Consider the restriction $A$ of $A_{0}\adjun$ to a subspace $\mathcal{D}$ containing $\ker B_{1} \cap \ker B_{2}$. Define a subspace of $\mathcal{B}^{2}$ by $\mathcal{C} := \bvec{B_{1} \\ B_{2}} \mathcal{D}$. Then the following claims are true
  \begin{enumerate}[label = \textrm{\textup{(\roman*)}}]
  \item The domain of $A$ can be written as
    \begin{align*}
      \dom A = \mathcal{D} = \set*{d \in \dom A_{0}\adjun : \bvec{B_{1} \\ B_{2}} d \in \mathcal{C}}
      .
    \end{align*}
  \item The operator closure of $A$ is $A_{0}\adjun$ restricted to
    \begin{align*}
      \mathcal{\tilde{D}} := \set*{d \in \dom A_{0}\adjun : \bvec{B_{1} \\ B_{2}} d \in \cl{\mathcal{C}}}
      ,
    \end{align*}
    where $\cl{\mathcal{C}}$ is the closure in $\mathcal{B}^{2}$. Therefore, $A$ is closed if and only if $\mathcal{C}$ is closed.
  \item The adjoint $A\adjun$ is the restriction of $-A_{0}\adjun$ to $\mathcal{D}'$, where
    \begin{align*}
      \mathcal{D}' := \set*{d' \in \dom A_{0}\adjun : \bvec{B_{1} \\ B_{2}} d' \in \bvec{0 & \idop \\ \idop & 0} \mathcal{C}^{\perp}}
      .
    \end{align*}
  \item The operator $A$ is (maximal) dissipative if and only if $\mathcal{C}$ is a (maximal) dissipative relation.    
  \end{enumerate}
\end{Proposition}

\section{Differential Operators}

Before we start analyzing port-Hamiltonian systems we will make some observation about the differential operators that will appear in the PDE. 
In this section we take care of all the technical details of these differential operators.
Since it doesn't really make a difference whether we use the scalar field $\R$ or $\C$ we will use $\K \in \set{\R,\C}$ for the scalar field.
The following assumption will be made for the rest of this work.

\begin{Voraussetzung}\label{assump:general-setting}
  Let $\irow,\icol,n\in\N$, $\Omega\subseteq \R^{n}$ be open with a bounded Lip\-schitz boundary, and $\myL = (L_{i})_{i=1}^{n}$ such that $L_{i}\in\K^{\irow\times \icol}$ for all $i \in \Zint{1}{n}$. Corresponding to $\myL$ we also have $\myL\hermitian := (L_{i}\hermitian)_{i=1}^{n}$, where $L_{i}\hermitian$ denotes the complex conjugated transposed (Hermitian transposed) matrix.
\end{Voraussetzung}

We will write $\mc D(\Omega)$ for the set of all $C^\infty(\Omega)$ functions with compact support in $\Omega$. Its dual space, the space of distribution, will be denoted by $\mc D\dual (\Omega)$. Moreover, we will write $\mc D(\R^n)\big\vert_{\Omega}$ for $\set{f\big\vert_{\Omega} : f \in \mc D(\R^n)}$. We will use $\partial_{i}$ as a short notation for $\frac{\partial}{\partial \zeta_{i}}$.

Sometimes it can be confusing to pay attention to the antilinear structure of an inner product of a Hilbert space, when switching between the inner product and the dual pairing. Thus, for the sake of clarity we will always consider the antidual space instead of the dual space, which is the space of all continuous antilinear mappings from the topological vector space into its scalar field. Hence, both the inner product and the (anti)dual pairing is linear in one component and antilinear in the other. So also $\mc D\dual(\Omega)$ is actually the antidual space of $\mc D(\Omega)$.

\begin{Definition}\label{def:diff-operator}
  Let $\myL$ be as in \autoref{assump:general-setting}. Then we define
  \begin{align*}
    \diffop := \sum_{i=1}^{n} \partial_{i} L_{i}
    \quad \text{and} \quad
    \diffopad := \sum_{i=1}^{n} \partial_{i} L\hermitian_{i}
  \end{align*}
  as operators on $\mathcal{D}\dual(\Omega)^{\icol}$ and $\mathcal{D}\dual(\Omega)^{\irow}$, respectively. Furthermore, we define the spaces
  \begin{align*}
    \Lspace &:= \set[\big]{ f \in \Lp{2}(\Omega,\K^{\icol}) : \diffop f \in \Lp{2}(\Omega, \K^{\irow}) } \\
    \mspace{-10mu}\text{and}\quad
    \Ladspace &:= \set[\big]{ f \in \Lp{2}(\Omega,\K^{\irow}) : \diffopad f \in \Lp{2}(\Omega, \K^{\icol}) }
                .
  \end{align*}
  These spaces are endowed with the inner product
  \begin{align*}
    \scprod{f,g}_{\Lspace} &:= \scprod{f,g}_{\Lp{2}(\Omega,\K^{\icol})} + \scprod{\diffop f, \diffop g}_{\Lp{2}(\Omega,\K^{\irow})} \\
    \mspace{-10mu}\text{and}\quad
    \scprod{f,g}_{\Ladspace} &:= \scprod{f,g}_{\Lp{2}(\Omega,\K^{\irow})} + \scprod{\diffopad f, \diffopad g}_{\Lp{2}(\Omega,\K^{\icol})}
  \end{align*}
  respectively. The space $\Lspace[0]$ is defined as $\cl[\Lnorm{.}]{\mc D(\Omega)^{\icol}}$ and $\Ladspace[0]$ analogously.
  We denote the trace operator by $\boundtr : H^{1}(\Omega) \to L^{2}(\partial \Omega)$ and the outward pointing normed normal vector on $\partial \Omega$ by $\nu$. 
  We define
  \begin{align*}
    \Lnub := \sum_{i=1}^{n} \nu_{i} L_{i}
    ,\quad \text{and} \quad
    \Lnubad := \sum_{i=1}^{n} \nu_{i} L_{i}\hermitian
    .
  \end{align*}
\end{Definition}

\begin{Remark}
  Clearly, $ H^{1}(\Omega,\K^{\icol})\subseteq \Lspace$ and  $H^{1}(\Omega,\K^{\irow})\subseteq \Ladspace$. It is also easy to see that $-\diffopad$ is the formal adjoint of $\diffop$.
\end{Remark}

For convenience we will write $H^{1}(\Omega)^{k}$ instead of $H^{1}(\Omega,\K^{k})$ and $\Lp{2}(\Omega)^{k}$ instead of $\Lp{2}(\Omega,\K^{k})$ for $k \in \N$.

\begin{Lemma}\label{le:closed-operator}
  The operator $\diffop$ with $\dom \diffop = \Lspace$ is a closed operator from $\Lp{2}(\Omega)^{\icol}$ to $\Lp{2}(\Omega)^{\irow}$ and $\Lspace$ endowed with the inner product $\scprod{.,.}_{\Lspace}$ is a Hilbert space.
\end{Lemma}

\begin{Proof}
  Let $\big((f_{n},\diffop f_{n})\big)_{n\in\N}$ be a sequence in $\diffop$ that converges to a point $(f,g)\in \Lp{2}(\Omega)^{\icol} \times \Lp{2}(\Omega)^{\irow}$.
  \newcommand{\scduala}[2][]{\dualprod[#1]{#2}_{\mathcal{D}\dual(\Omega)^{\icol},\mathcal{D}(\Omega)^{\icol}}}%
  \newcommand{\scdualb}[2][]{\dualprod[#1]{#2}_{\mathcal{D}\dual(\Omega)^{\irow},\mathcal{D}(\Omega)^{\irow}}}%
  For an abitrary $\phi \in \mathcal{D}(\Omega)^{\irow}$ we have
  \begin{align*}
    \scdualb{g, \phi} & = \lim_{n\in\N} \scdualb{\diffop f_{n},\phi} \\
                      & = \lim_{n\in\N}\scduala{f_{n}, -\diffopad \phi} \\
                      & = \scduala{f, -\diffopad \phi} \\
                      & = \scdualb{\diffop f, \phi},
  \end{align*}
  which implies $g = \diffop f$. Since $g$ is also in $\Lp{2}(\Omega)^{\irow}$, we conclude that $\diffop$ is closed.
\end{Proof}

\begin{Example}\label{ex:div-grad}
  Let us regard the following matrices
\begin{align*}
  L_{1} = \begin{bmatrix} 1 & 0 & 0 \end{bmatrix}
  ,\quad
  L_{2} = \begin{bmatrix} 0 & 1 & 0 \end{bmatrix}
  ,\quad\text{and}\quad
  L_{3} = \begin{bmatrix} 0 & 0 & 1 \end{bmatrix}
  .
\end{align*}
Then we obtain the corresponding differential operators
\begin{align*}
  \diffop = \begin{bmatrix}  \partial_{1} & \partial_{2} & \partial_{3} \end{bmatrix} = \Div
  \quad \text{and}\quad
  \diffopad = \begin{bmatrix}  \partial_{1} \\ \partial_{2} \\ \partial_{3} \end{bmatrix} = \grad
  .
\end{align*}
The corresponding operator $\Lnub$ that acts on $\Lp{2}(\partial\Omega)$ can be written as an inner product
\begin{equation*}
  \Lnub f
  =
  \begin{bmatrix}
    \nu_1 & \nu_2 & \nu_3
  \end{bmatrix}
  \begin{bmatrix}
    f_1 \\ f_2 \\ f_3
  \end{bmatrix}
  =
  \nu \cdot f = \scprod{f,\nu}_{\K^3}
  .
\end{equation*}
\end{Example}

Clearly the previous example can be extended to any finite dimension.

\begin{Example}\label{ex:rot}
  The following matrices will construct the rotation operator.
  \begin{align*}
    L_{1}=
    \begin{bmatrix}
      0 & 0 & 0 \\ 0 & 0 & -1 \\ 0 & 1 & 0
    \end{bmatrix}
    ,\quad
    L_{2}=
    \begin{bmatrix}
      0 & 0 & 1 \\ 0 & 0 & 0 \\ -1 & 0 & 0
    \end{bmatrix}
    ,\quad\text{and}\quad
    L_{3}=
    \begin{bmatrix}
      0 & -1 & 0 \\ 1 &0 &0 \\ 0&0&0
    \end{bmatrix}.
  \end{align*}
  In this example we have $L_{i}\hermitian = - L_{i}$. Furthermore, the corresponding differential operator is
  \begin{align*}
    \diffop =
    \begin{bmatrix}
      0 & -\partial_{3} & \partial_{2} \\
      \partial_{3} & 0 & -\partial_{1} \\
      -\partial_{2} & \partial_{1} & 0
    \end{bmatrix}
    = \rot = -\diffopad
    .
  \end{align*}
  The corresponding operator $\Lnub$ that acts on $\Lp{2}(\partial\Omega)$ can be written as a cross product
  \begin{equation*}
    \Lnub f
    =
    \begin{bmatrix}
      0 & -\nu_{3} & \nu_{2} \\
      \nu_{3} & 0 & -\nu_{1} \\
      -\nu_{2} & \nu_{1} & 0
    \end{bmatrix}
    \begin{bmatrix}
      f_1 \\ f_2 \\ f_3
    \end{bmatrix}
    = \nu \times f
    .
  \end{equation*}
\end{Example}

\begin{Lemma}
  The adjoint of $\diffop$ with $\dom \diffop = \Lspace$ is $\diffop\adjun = -\diffopad$ on $\dom \diffop\adjun \subseteq \Ladspace$.
\end{Lemma}

\begin{Proof}
   \newcommand{\scdual}[2][]{\dualprod[#1]{#2}_{\mathcal{D}\dual,\mathcal{D}}}%
   For an arbitrary $g \in \dom \diffop\adjun$ and an arbitrary $\phi \in \mathcal{D}(\Omega)$ we have (to shorten the notation we will write $\scdual{.,.}$ instead of $\dualprod{.,.}_{\mathcal{D}\dual(\Omega)^{k},\mathcal{D}(\Omega)^{k}}$ and $\scprod{.,.}_{\Lp{2}}$ instead of $\scprod{.,.}_{\Lp{2}(\Omega)^{k}}$)
  \begin{align*}
    \scdual{\diffop\adjun g, \phi}
    = \scprod{\diffop\adjun g, \phi}_{\Lp{2}} = \scprod{g, \diffop \phi}_{\Lp{2}}
      = \scdual{g, \diffop \phi}
     = \scdual{-{ \diffopad g}, \phi}
      .
  \end{align*}
  Therefore, $\diffop\adjun g = - \diffopad g$ and $\diffop\adjun g \in \Lp{2}(\Omega)$ implies $\diffopad g \in \Lp{2}(\Omega)$. Consequently, $\dom \diffop\adjun \subseteq \Ladspace$. 
\end{Proof}

\begin{Remark}\label{re:skew-symmetric}
  If $\myL$ contains only Hermitian matrices ($L_{i}\hermitian = L_{i}$), then $\diffop\adjun$ is skew-symmetric.
\end{Remark}

\begin{Remark}\label{re:Lspace-convergence}
  The mapping $\iota: H(\diffop, \R^n) \to \Lspace, f \mapsto f\big\vert_{\Omega}$ is well-defined and continuous for any open set $\Omega\subseteq \R^n$. In particular, $\diffop(f\big\vert_{\Omega}) = (\diffop f)\big\vert_{\Omega}$. Hence, we can always regard an $f \in H(\diffop, \R^n)$ as an element of $\Lspace$, especially when $\supp f \subseteq \cl{\Omega}$. Moreover, if $f_n \to f$ in $H(\diffop, \R^n)$, then $f_n \to f$ in $\Lspace$.
\end{Remark}

\begin{Definition}
  A set $O\subseteq \R^{n}$ is \emph{strongly star-shaped} with respect to $x_{0}$, if for every $x\in\cl{O}$ the half-open line segment $\set{x_{0}+\theta x: \theta \in [0,1)}$ is contained in $O$. We call $O$ strongly star-shaped, if there is a $x_{0}$ such that $O$ is strongly star-shaped with respect to $x_{0}$.
\end{Definition}


Note that this is equivalent to
\[
  \theta(\cl{O} - x_{0}) + x_{0} \subseteq O \quad \text{for all}\quad \theta \in [0,1).
\]

\begin{Lemma}\label{le:star-shaped}
  Let $f\in H(\diffop, \R^n)$ and $x_0\in \R^n$. Furthermore, let $f_{\theta}(x) := f(\frac{1}{\theta}(x-x_{0}) + x_{0})$ for $\theta \in (0,1)$. Then $f_{\theta} \in H(\diffop, \R^n)$ and $f_{\theta} \to f$ in $H(\diffop, \R^n)$ as $\theta \to 1$.
  If there exisits a strongly star-shaped set $O$ such that $\supp f \subseteq \cl{O}$, then $\supp f_{\theta} \subseteq O$ for $\theta \in (0,1)$.
\end{Lemma}

\begin{Proof}
  Let $\alpha(x) := \frac{1}{\theta} (x-x_{0}) + x_{0}$. Then we have $f_{\theta} = f\circ \alpha$ and
    \begin{align*}
    \dualprod{\diffop (f\circ \alpha), \phi}_{\mc D\dual, \mc D} &= \dualprod{f, -(\diffopad \phi)\circ \alpha^{-1} \theta^n}_{\mc D\dual, \mc D}
    = \dualprod[\Big]{f, -\sum_{i=1}^n L_i\hermitian \partial_i \Big(\phi\circ \alpha^{-1}\frac{1}{\theta}\Big) \theta^n}_{\mc D\dual, \mc D} \\
    &= \dualprod[\Big]{f, -\diffopad \Big(\frac{1}{\theta} \phi\circ\alpha^{-1}\Big) \theta^n}_{\mc D\dual, \mc D}
      = \dualprod[\Big]{\frac{1}{\theta} (\diffop f)\circ\alpha, \phi}_{\mc D\dual, \mc D}
      .
    \end{align*}
    Therefore, $\diffop f_{\theta} = \frac{1}{\theta} (\diffop f)_{\theta}$ and $f_{\theta} \in H(\diffop, \R^n)$. We can also write $f_{\theta}$ as \linebreak $T_{x_{0}} D_{\frac{1}{\theta}} T_{-x_{0}} f$ where $T_{y}$ is the translation mapping $f \mapsto f(.+y)$ and $D_{\eta}$ is the dilation mapping $f \mapsto f(\eta.)$. Since $T_{y}$ is bounded and $D_{\eta}$ converges strongly to $\idop$ as $\eta \to 1$, we conclude $f_{\theta} \to f$ in $\Lp{2}(\R^n)^{\icol}$ as $\theta \to 1$ and $\diffop f_{\theta} = \frac{1}{\theta} (\diffop f)_{\theta} \to \diffop f$ in $\Lp{2}(\R^{n})^{\irow}$ as $\theta \to 1$. Hence, $f_{\theta} \to f$ in $H(\diffop, \R^n)$.

  Let $O$ be strongly star-shaped with respect to $x_0$ and $\supp f \subseteq \cl{O}$. Then for $\theta \in (0,1)$
  \begin{equation*}
    \supp f_{\theta} = \theta (\supp f - x_0) + x_0 \subseteq \theta (\cl{O} - x_0) + x_0 \subseteq O
    .
    \qedhere
  \end{equation*}
\end{Proof}

\begin{Remark}
  If $f \in \Lspace$ and $\psi \in \mc D(\R^n)\big\vert_{\Omega}$, then by the product rule for distributional derivatives also $\psi f \in \Lspace$ and $\diffop (\psi f) = \psi \diffop f + \sum_{i=1}^n (\partial_i \psi) L_i f$.
\end{Remark}

\begin{Lemma}\label{th:approx-with-compact-supp}
  For every $f\in H(\diffop, \R^n)$ exists a sequence $(f_k)_{k\in\N}$ in $H(\diffop, \R^n)$ with $\supp f_k$ is compact that converges to $f$ in $H(\diffop, \R^n)$.
\end{Lemma}

\begin{Proof}
  Let $\psi \in C^\infty(\R^n,\R)$ be such that
  \[
    \psi(\zeta) =
    \begin{cases}
      1, & \text{if}\; \norm{\zeta} \leq 1, \\
      0, & \text{if}\; \norm{\zeta} \geq 2.
    \end{cases}
  \]
  Then $f_k := \psi(\frac{1}{k} .) f \in H(\diffop, \R^n)$ and $f_k \to f$ in $\Lp{2}$. By $\diffop f_k = \psi(\frac{1}{k}.) \diffop f + \frac{1}{k}\sum_{i=1}^n(\partial_i\phi)(\frac{1}{k}.)L_if$, we conclude $f_k \to f$ in $H(\diffop, \R^n)$.
\end{Proof}

The next lemma is similiar to \cite[Lemma 1, page 206]{book-DL} the main idea of the proof can be adopted.

\begin{Lemma}\label{le:Lspace0}
  If $f \in \Lspace$ is such that
  \begin{align}\label{eq:zero-boundary-condition}
    \scprod{\diffop f, \phi}_{\Lp{2}(\Omega)} + \scprod{f, \diffopad \phi}_{\Lp{2}(\Omega)} = 0
    \quad \text{for all}\quad
    \phi \in \mathcal{D}(\R^{n})^{\irow},
  \end{align}
  then $f \in \Lspace[0]$.
\end{Lemma}

\begin{Proof}
  Let $f \in \Lspace$ such that it satisfies \eqref{eq:zero-boundary-condition}. Then we have to find a sequence $(f_n)_{n\in\N}$ in $\mc D(\Omega)^{\icol}$ that converges to $f$ with respect to $\Lnorm{.}$.

  We define $f_0 := \diffop f$ and $\tilde{f}, \tilde{f}_0$ as the extension of $f$ and $f_0$ respectively on $\R^n$ such that these functions are $0$ outside of $\Omega$. By
  \begin{align*}
    \dualprod{\tilde{f}_0, \phi}_{\mc D\dual(\R^n), \mc D(\R^n)}
    &= \scprod{\tilde{f}_{0}, \phi}_{\Lp{2}(\R^n)} = \scprod{\diffop f, \phi}_{\Lp{2}(\Omega)} \stackrel{\eqref{eq:zero-boundary-condition}}{=}  \scprod{f,-\diffopad\phi}_{\Lp{2}(\Omega)} \\
    &= \scprod{\tilde{f}, -\diffopad\phi}_{\Lp{2}(\R^n)}
    =  \dualprod{\tilde{f}, -\diffopad\phi}_{\mc D\dual(\R^n), \mc D(\R^n)}
  \end{align*}
  for $\phi \in \mc D(\R^n)^{\irow}$, we see that $\tilde{f}_0 = \diffop \tilde{f}$ and $\tilde{f}\in H(\diffop,\R^n)$ with $\supp \tilde{f} \subseteq \cl{\Omega}$.

  \medskip\par\noindent
  \emph{Step 1. Assume that there is a bounded $\Omega'\subseteq \Omega$ with bounded Lipschitz boundary such that $\supp \tilde{f} \subseteq \cl{\Omega'}$.} By \cite[Proposition 2.5.4, page 69]{star-shaped-covering} there is a finite open covering $(O_i)_{i=1}^{k}$ of $\cl{\Omega'}$ such that $O_{i} \cap \Omega'$ is strongly star-shaped. We employ a partition of unity and obtain $(\alpha_i)_{i=1}^{k}$, subordinate to this covering, that is
  \begin{equation*}
    \alpha_i \in \mc D(O_i), \quad \alpha_i(x) \in [0,1] \quad\text{and}\quad \sum_{i=1}^{k} \alpha_i(x) = 1
    \quad\text{for all}\quad x\in \Omega'
    .
  \end{equation*}
  Hence, $\tilde{f} =  \sum_{i=1}^{k} \alpha_i \tilde{f}$ and we define $f_i := \alpha_{i}\tilde{f}$. By construction $f_i \in H(\diffop, \R^n)$ and $\supp f_i \subseteq \cl{O_i \cap \Omega'}$, where $O_i \cap \Omega'$ is strongly star-shaped. \autoref{le:star-shaped} ensures that $\supp (f_i)_\theta \subseteq O_i \cap \Omega'$ for $\theta \in (0,1)$ and $(f_i)_\theta \to f_i$ in $H(\diffop,\R^n)$ for $\theta \to 1-$.
  
  Let $\rho_{\epsilon}$ be a positive $C^{\infty}$ mollifier with compact support. Then $\rho_{\epsilon} \ast g \to g$ in $\Lp{2}(\R^n)$ for an arbitrary $g \in \Lp{2}(\R^n)$. Since $\diffop (\rho_{\epsilon} \ast h) = \rho_{\epsilon} \ast \diffop h$, we also have that $\rho_{\epsilon} \ast h \to h$ in $H(\diffop, \R^n)$ for $h\in H(\diffop, \R^n)$ and $\rho_{\epsilon} \ast h \in \mc D(\R^n)^{\icol}$.

  For fixed $\theta \in (0,1)$ and $\epsilon$ sufficiently small, we can say $\supp \rho_{\epsilon}\ast (f_i)_\theta \subseteq \Omega'$. This establishes the existence of a sequence $(\rho_{\epsilon_j} \ast (f_{i})_{\theta_j})_{j\in\N}$ in $\mc D(\Omega')^{\icol}$ converging to $f_i$ in $H(\diffop, \R^n)$.
  Doing this for every $i \in \Zint{1}{k}$ yields sequences $(f_{i,j})_{j\in\N}$ in $\mc D(\Omega')^{\icol}\subseteq \mc D(\Omega)^{\icol}$ converging to $f_i$ in $H(\diffop, \R^n)$. Consequently, $\big(\sum_{i=1}^{k} f_{i,j}\big)_{j\in\N}$ is a sequence in $\mc D(\Omega)^{\icol}$ that converges to $\tilde{f}$ in $H(\diffop, \R^n)$ and by \autoref{re:Lspace-convergence} also in $\Lspace$.

  \medskip\par\noindent
  \emph{Step 2. Without extra assumptions.} By the already shown each entry of the sequence $(f_k)_{k\in\N}$ from \autoref{th:approx-with-compact-supp} can be approximated by $\mc D(\Omega)^{\icol}$ elements. A diagonalization argument yields the same for the limit $\tilde{f}$. By \autoref{re:Lspace-convergence} this diagonal sequence also converges in $\Lspace$.
\end{Proof}

\begin{Theorem}
  $\mc D(\R^{n})^{\icol}\big\vert_{\Omega}$ is dense in $\Lspace$.
\end{Theorem}

\begin{Proof}
  Suppose $\mc D(\R^{n})^{\icol}\big\vert_{\Omega}$ is not dense in $\Lspace$. Then there exists a non zero $f\in\Lspace$ such that
  \begin{align}\label{eq:ortho-on-dense}
    \scprod{f, g}_{\Lspace} = \scprod{f,g}_{\Lp{2}} + \scprod{\diffop f, \diffop g}_{\Lp{2}} = 0
    \quad \text{for all} \quad g \in \mc D(\R^{n})^{\icol}\big\vert_{\Omega}
    .
  \end{align}
  In particular, for an arbitrary $h \in \mathcal{D}({\Omega})^{\icol}$ we have
  \begin{align*}
    \dualprod{f,h}_{\mathcal{D}\dual, \mathcal{D}} = \scprod{f,h}_{\Lp{2}}
    = - \scprod{\diffop f, \diffop h}_{\Lp{2}}
    = - \dualprod{\diffop f, \diffop h}_{\mathcal{D}\dual, \mathcal{D}}
    = \dualprod{\diffopad \diffop f, h}_{\mathcal{D}\dual, \mathcal{D}},
  \end{align*}
  which implies that $f = \diffopad\diffop f \in \Lp{2}(\Omega)^{\icol}$. Hence we can rewrite \eqref{eq:ortho-on-dense} as
  \begin{align*}
    \scprod{\diffopad \underbrace{\diffop f}_{=f_{0}},g}_{\Lp{2}(\Omega)} + \scprod{\underbrace{\diffop f}_{=f_{0}}, \diffop g}_{\Lp{2}(\Omega)} = 0
    \quad \text{for all} \quad g \in \mc D(\R^{n})^{\icol}\big\vert_{\Omega}
    .
  \end{align*}
  By \autoref{le:Lspace0} (switching the roles of $\diffop$ and $\diffopad$), defining $f_{0} := \diffop f$ yields $f_{0} \in \Ladspace[0]$. Since $\mathcal{D}({\Omega})^{\irow}$ is dense in $\Ladspace[0]$, there is a sequence $(f_{n})_{n\in\N}$ in $\mathcal{D}({\Omega})^{\irow}$ converging to $f_{0}$ with respect to $\norm{.}_{\Ladspace}$. Note that $f = \diffopad f_{0}$. 
  \begin{align*}
    \scprod{f_{0}, f_{n}}_{\Ladspace}
    & = \scprod{f_{0}, f_{n}}_{\Lp{2}} + \scprod{\diffopad f_{0}, \diffopad f_{n}}_{\Lp{2}}
      = \scprod{\diffop f, f_{n}}_{\Lp{2}} + \scprod{f , \diffopad f_{n}}_{\Lp{2}} \\
    & = \dualprod{\diffop f, f_{n}}_{\mathcal{D}\dual,\mathcal{D}}
      - \dualprod{\diffop f, f_{n}}_{\mathcal{D}\dual,\mathcal{D}} = 0.
  \end{align*}
  Since $\norm{f_{0}}_{\Ladspace}^{2} = \lim_{n\in\N} \scprod{f_{0}, f_{n}}_{\Ladspace} = 0$, we have that $f_{0} = 0$, which implies $f = \diffopad f_0 = 0$. Hence, $\mathcal{D}(\R^{n})^{\icol}\big\vert_{\Omega}$ is dense in $\Lspace$.
\end{Proof}

\begin{Lemma}\label{le:int-by-parts}
 Let $f \in H^{1}(\Omega)^{\icol}$ and $g\in H^{1}(\Omega)^{\irow}$. Then we have
  \begin{align}
  \begin{split}\label{eq:int-by-parts}
    \scprod{\diffop f,g}_{\Lp{2}(\Omega)^{\irow}} + \scprod{f,\diffopad g}_{\Lp{2}(\Omega)^{\icol}}
    	&
	= \scprod{\Lnub \boundtr f,\boundtr g}_{\Lp{2}(\partial \Omega)^{\irow}} 
	\\
    	&= \scprod{\boundtr f,\Lnubad \boundtr g}_{\Lp{2}(\partial \Omega)^{\icol}}.
  \end{split}
  \end{align}
\end{Lemma}

\begin{Proof}
By the definition of $\diffop$ and $\diffopad$, and the linearity of the scalar product we can write the left-hand-side of \eqref{eq:int-by-parts} as
  \begin{align*}
    \int_{\Omega} \sum_{i=1}^{n} \scprod{\partial_{i}L_{i}f,g} + \scprod{f,\partial_{i} L_{i}\hermitian g}\dx[\la]
    = \int_{\Omega} \sum_{i=1}^{n} \scprod{\partial_{i}L_{i}f,g} + \scprod{L_{i} f,\partial_{i} g}\dx[\la]
    ,
  \end{align*}
where $\la$ denotes the Lesbesgue measure.
By the product rule for derivatives and the Gau{\ss}'s theorem (divergence theorem) this is equal to
  \begin{align*}
    \int_{\Omega} \sum_{i=1}^{n} \partial_{i} \scprod{L_{i}f,g} \dx[\la]
    = \int_{\partial \Omega} \sum_{i=1}^{n} \nu_{i} \boundtr \scprod{ L_{i}f,  g} \dx[\mu]
    = \int_{\partial \Omega}  \scprod{\Lnub \boundtr f, \boundtr g} \dx[\mu]
    ,
  \end{align*}
  where $\nu$ denotes the  outward pointing normed normal vector on $\partial \Omega$ and $\mu$ denotes the surface measure of $\partial \Omega$.
\end{Proof}

\begin{Corollary}\label{co:int-by-parts}
  Let $f \in H^{1}(\Omega)^{\icol}$ and $g\in H^{1}(\Omega)^{\irow}$. Then we have
  \begin{align*}
    \abs[\Big]{\scprod{\Lnub \boundtr f,\boundtr g}_{\Lp{2}(\partial \Omega)^{\irow}}}
    \leq
    \Lnorm{f} \Ladnorm{g}
    .
  \end{align*}
\end{Corollary}

\begin{Proof}
  \autoref{le:int-by-parts}, the triangular inequality and Cauchy Schwartz's inequality yield
  \begin{align*}
    \abs[\Big]{\scprod{\Lnub \boundtr f,\boundtr g}_{\Lp{2}(\partial \Omega)^{\irow}}}
    &\leq
      \abs[\big]{\scprod{\diffop f,g}_{\Lp{2}(\Omega)^{\irow}}} + \abs[\big]{\scprod{f,\diffopad g}_{\Lp{2}(\Omega)^{\icol}}}
    \\
    &\leq
      \norm{\diffop f} \norm{g} + \norm{f} \norm{\diffopad g}
    \\
    &\leq
      \sqrt{\norm{\diffop f}^{2} + \norm{f}^{2}} \sqrt{\norm{g}^{2} + \norm{\diffopad g}^{2}} 
    \\
    &= \Lnorm{f}\Ladnorm{g}
      .
      \qedhere
  \end{align*}
\end{Proof}

\section{Quasi Gelfand Triples}\label{sec:quasi-gelfand-triple}

Normally when we talk about Gelfand triples we have a Hilbert space $H_{0}$ and another Hilbert space $H$ that can be continuously and densely embedded into $H_{0}$. We want to weaken this condition such that the norm of $H$ isn't necessarily related to the norm of $H_{0}$.
\medskip\par
We will have the following setting: Let $(\hsmid,\scprod{.,.}_{\hsmid})$ be a Hilbert space and $\scprod{.,.}_{\hs}$ another inner product (not necessarily related to $\scprod{.,.}_{\hsmid}$) which is defined on a dense (w.r.t. $\norm{.}_{\hsmid}$) subspace $\tilde{D}_{+}$ of $\hsmid$.

We denote the completion of $\tilde{D}_{+}$ w.r.t. $\norm{.}_{\hs} = \sqrt{\scprod{.,.}_{\hs}}$ by $\hs$. This completion is again a Hilbert space with the extension of $\scprod{.,.}_{\hs}$, for which we use the same symbol. Now we have that $\tilde{D}_{+}$ is dense in $\hsmid$ w.r.t. $\norm{.}_{\hsmid}$ and dense in $\hs$ w.r.t. $\norm{.}_{\hs}$.

\begin{Definition}\label{def:Dminus}
  Let $\hsmid$, $\hs$ and $\tilde{D}_{+}$ be as mentioned in the beginning of this section. Then we define
  \begin{align*}
    \norm{g}_{\hsdual} := 
    \sup_{f\in \tilde{D}_{+}\setminus\set{0}} \frac{%
    \abs{\scprod{g,f}_{\hsmid}}
    }{\norm{f}_{\hs}}
    \quad \text{and} \quad
    D_{-} := \set[\Big]{g \in \hsmid : \norm{g}_{\hsdual} < +\infty}.
  \end{align*}
  We denote the completion of $D_{-}$ w.r.t. $\norm{.}_{\hsdual}$ by $\hsdual$.
\end{Definition}

\begin{Remark}
By definition of $D_{-}$ we can identify every $g \in D_{-}$ with an element of $\hs\dual$ by the continuous extension of $ f \in \tilde{D}_{+} \mapsto \scprod{g,f}_{\hsmid}$ to $\hs$. The completion $\hsdual$ is isomorphic to the closure of $D_{-}$ in $\hs\dual$.
\end{Remark}

\begin{Lemma}\label{le:Dminus-complete}
  $D_{-}$ is complete with respect to $\norm{g}_{\hsdual \cap \hsmid}^{2} := \norm{g}_{\hsmid}^{2} + \norm{g}_{\hsdual}^{2}$.
\end{Lemma}

\begin{Proof}
  Let $(g_{n})_{n\in\N}$ be a Cauchy sequence in $D_{-}$ with respect to $\norm{.}_{\hsdual\cap\hsmid}$. Then $(g_{n})_{n\in\N}$ is a convergent sequence in $\hsmid$ (w.r.t. $\norm{.}_{\hsmid}$) and a Cauchy sequence in $D_{-}$ (w.r.t. $\norm{.}_{\hsdual}$). We denote the limit in $\hsmid$ by $g_{0}$. We obtain
  \begin{align*}
    \abs{\scprod{g_{0},f}_{\hsmid}} = \lim_{n\in\N} \abs{\scprod{g_{n},f}_{\hsmid}}
    \leq \lim_{n\in\N} \norm{g_{n}}_{\hsdual} \norm{f}_{\hs} \leq C \norm{f}_{\hs}
  \end{align*}
  and consequently $g_{0} \in D_{-}$.

  Let $\epsilon > 0$ be arbitray. Since $(g_n)_{n\in\N}$ is a Cauchy sequence with respect to $\norm{.}_{\hsdual}$, there is an $n_{0}\in \N$ such that for all $f\in \tilde{D}_{+}$ with $\norm{f}_{\hs} = 1$
  \begin{align*}
    \abs{\scprod{g_{n}-g_{m}, f}} \leq \frac{\epsilon}{2}, \quad \text{if}\quad n,m \geq n_{0}
  \end{align*}
  holds true. Furthermore, for all $f \in \tilde{D}_{+}$ there exists an $m_{f}\geq n_{0}$ such that $\abs{\scprod{g_0 - g_{m_{f}},f}} \leq \frac{\epsilon \norm{f}_{\hs}}{2}$. This yields
  \begin{align*}
    \frac{\abs{\scprod{g_{0}-g_{n},f}}}{\norm{f}_{\hs}} \leq \frac{\abs{\scprod{g_{0} - g_{m_{f}},f}}}{\norm{f}_{\hs}} + \frac{\abs{\scprod{g_{m_{f}}-g_{n},f}}}{\norm{f}_{\hs}} \leq \epsilon.
  \end{align*}
  Since the right-hand-side is independent of $f$, we obtain
  \begin{align*}
   \norm{g_{0}-g_{n}}_{\hsdual} = \sup_{f\in \tilde{D}_{+}} \frac{\abs{\scprod{g_{0}-g_{n},f}}}{\norm{f}_{\hs}} \leq \epsilon.
  \end{align*}
  Hence, $g_{0}$ is also the limit of $(g_{n})_{n\in\N}$ with respect to $\norm{.}_{\hsdual}$ and consequently the limit of  $(g_{n})_{n\in\N}$ with respect to $\norm{.}_{\hsdual\cap\hsmid}$.
\end{Proof}

\begin{Lemma}
  The embedding $\iota_{+}: \tilde{D}_{+} \subseteq \hs \to \hsmid, f \mapsto f$ is a dense defined operator with $\ran \iota_{+}$ is dense and $\ker \iota_{+} = \set{0}$.
\end{Lemma}

\begin{Proof}
  By assumption on $\tilde{D}_{+}$ the embedding $\iota_{+}$ is dense defined and has a dense range. Clearly, $\ker \iota_{+} = \set{0}$.
\end{Proof}

\begin{Lemma}\label{le:adjoint-embedding}
  Let $\iota_{+}\adjun$ denote the adjoint relation of the embedding mapping $\iota_{+}$ in the previous lemma.
  Then $\iota_{+}\adjun$ is single-valued ($\mul \iota_{+}\adjun = \set{0}$) and $\ker \iota_{+}\adjun = \set{0}$. Its domain conincides with $D_{-}$ and $\iota_{+}\adjun : D_{-} \to \hs$ is isometric w.r.t. $\norm{.}_{\hsdual}$.
  
  If $\ker \cl{\iota_{+}} = \set{0}$, then $\ran \iota\adjun$ is dense in $\hs$.
\end{Lemma}

\begin{Proof}
  The density of the domain of $\iota_{+}$ yields $\mul \iota_{+}\adjun = (\dom \iota_{+})^{\perp} = \set{0}$, and $\cl[\hsmid]{\ran \iota_{+}} = \hsmid$ yields $\ker \iota_{+}\adjun = \set{0}$.
  The following equivalences show that $\dom \iota_{+} = D_{-}$:
\begin{align*}
  g \in \dom \iota_{+}\adjun &\Leftrightarrow \scprod{g,\iota_{+} f}_{\hsmid} \; \text{is continuous in} \; f \in \tilde{D}_{+} \;\text{w.r.t.}\; \norm{.}_{\hs} \\
                         &\Leftrightarrow \sup_{f\in \tilde{D}_{+}} \frac{\abs{\scprod{g, f}_{\hsmid}}}{\norm{f}_{\hs}} < +\infty \\
                         &\Leftrightarrow g \in D_{-}.
\end{align*}
For $g\in D_{-}$ we have
\begin{align*}
  \norm{g}_{\hsdual} = \sup_{f\in \tilde{D}_{+}} \frac{\abs{\scprod{g,f}_{\hsmid}}}{\norm{f}_{\hs}} = \sup_{f\in \tilde{D}_{+}} \frac{\abs{\scprod{\iota_{+}\adjun g,f}_{\hs}}}{\norm{f}_{\hs}} = \norm{\iota_{+}\adjun g}_{\hs},
\end{align*}
which proves that $\iota_{+}\adjun$ is isometric.

If $\ker \cl{\iota_{+}} = \set{0}$, then the following equation yields the density of $\ran \iota_{+}\adjun$ in $\hs$
\begin{equation*}
\ker \cl{\iota_{+}} = \ker \iota_{+}\adjun[2] = (\ran \iota_{+}\adjun)^{\perp}.
\qedhere
\end{equation*}
\end{Proof}

\begin{Proposition}\label{prop:iota-closed-equivalences}
  The following assertions are equivalent.
  \begin{enumerate}[label = \textrm{\textup{(\roman*)}}]
    \item There is a topological vector space $(Z,\mathcal{T})$ that contains $\hsmid$ and $\hs$ such that $\tilde{D}_+ \subseteq \hs\cap \hsmid$ in $Z$, and the topology $\mathcal{T}$ is coarser (weaker) than the topology of $\norm{.}_{\hsmid}$ and coarser (weaker) than the topology of $\norm{.}_{\hs}$.
    \item If $\tilde{D}_{+}\ni f_n \to 0$ w.r.t. $\norm{.}_{\hs}$ and $\lim_{n\in\N} f_n$ exists w.r.t. $\norm{.}_{\hsmid}$ then this limit is also $0$ and if $\tilde{D}_{+}\ni f_n \to 0$ w.r.t. $\norm{.}_{\hsmid}$ and $\lim_{n\in\N} f_n$ exists w.r.t. $\norm{.}_{\hs}$ then this limit is also $0$.
    \item $\iota_{+}: \tilde{D}_{+} \subseteq \hs \to \hsmid, f \mapsto f$ is closable and its closure is injective.
    \item $D_{-}$ is dense in $\hsmid$ and dense in $\hs\dual$.
  \end{enumerate}
\end{Proposition}

\begin{Proof}
  $\textup{(i)} \Rightarrow \textup{(ii)}$:
  Let $(f_n)_{n\in\N}$ be a sequence in $\tilde{D}_{+}$ such that $f_n \to 0$ w.r.t. $\hs$ and $f_n \to f$ w.r.t. $\hsmid$. Since $\mc T$ is coarser than both of the topologies induced by these norms, we also have
  \[
  \begin{tikzcd}[row sep=0.1em]
    &  0 \\
    f_n \arrow{ur}{\mc T}\arrow{dr}{\mc T} & \\
    & f 
  \end{tikzcd}
  \quad \text{in } Z.
  \]
  Since $\mc T$ is Hausdorff, we conclude $f=0$. Analogously, we can show the converse statement.
  
  \smallskip\par
  $\textup{(ii)} \Rightarrow \textup{(iii)}$: If $(f_n,f_n)_{n\in\N}$ is a sequence in $\iota_{+}$ that converges to $(0,f)\in\hs\times \hsmid$, then $f = 0$ by $(ii)$. Hence, $\mul \cl{\iota_{+}} = \set{0}$ and consequently $\iota_{+}$ is closable. Analogously, we can show that $\ker \cl{\iota_{+}} = \set{0}$.
  \smallskip\par
  $\textup{(iii)} \Rightarrow \textup{(iv)}$: We have $(\dom \iota_{+}\adjun)^{\perp} = \mul \cl{\iota_{+}}$, where $\cl{\iota_{+}}$ is the closure of $\iota_{+}$. Since $\iota_{+}$ is closable, we have $\mul \cl{\iota_{+}} = \set{0}$, which yields $\dom \iota_{+}\adjun$ is dense. By \autoref{le:adjoint-embedding} $\dom \iota_{+}\adjun$ coincides with $D_{-}$. The second assertion of \autoref{le:adjoint-embedding} yields that $D_{-}$ is dense in $\hs\dual$.
  \smallskip \par
  $\textup{(iv)} \Rightarrow \textup{(i)}$: By \autoref{le:Dminus-complete} $D_{-}$ is complete with respect to $\norm{g}_{\hsdual\cap\hsmid}^{2} := \norm{g}_{\hsdual}^{2} + \norm{g}_{\hsmid}^{2}$. Since $D_{-}$ dense in $\hsmid$ and embedding into $\hsmid$ is continuous, we can construct an ordinary Gelfand triple. Hence, $Z$, the completion of $\hsmid$ with respect to $\norm{z}_{Z} := \sup_{g\in D_{-}\setminus \set{0}} \frac{\abs{\scprod{z,g}_{\hsmid}}}{\norm{g}_{\hsdual\cap\hsmid}}$, is the dual space of $D_{-}$ with respect to the pivot space $\hsmid$. For $f\in \tilde{D}_{+}$ we have
  \begin{align*}
    \norm{f}_{Z} = \sup_{g\in D_{-}\setminus \set{0}} \frac{\abs{\scprod{f,g}_{\hsmid}}}{\norm{g}_{\hsdual\cap\hsmid}}
    \leq \sup_{g\in D_{-}\setminus \set{0}} \frac{\norm{f}_{\hs}\norm{g}_{\hsdual}}{\norm{g}_{\hsdual\cap\hsmid}}
    \leq \norm{f}_{\hs}.
  \end{align*}
  Consequently, we can regard $\hs$ as a subspace of $Z$. By contruction the topology of $\norm{.}_{Z}$ is coarser than the topology of $\norm{.}_{\hsmid}$ and by the last inequality it is also coarser than the topology of $\norm{.}_{\hs}$.
\end{Proof}

From now on we will assume that one and therefore all properties in \autoref{prop:iota-closed-equivalences} are satisfied. Therefore, $\hs \cap \hsmid$ is well-defined and complete with the norm $\norm{.}_{{\hs\cap\hsmid}} := \sqrt{\norm{.}_{\hs}^{2} + \norm{.}_{\hsmid}^{2}}$.

\begin{Lemma}
  $\tilde{D}_{+}$ is dense in $\hs\cap\hsmid$ with respect to $\norm{.}_{\hs\cap\hsmid}^{2}:= \norm{.}_{\hs}^{2} + \norm{.}_{\hsmid}^{2}$.
\end{Lemma}

\begin{Proof}
  We define $P_{+} := \hs \cap \hsmid$ and we define $P_{-}$ analogously to $D_{-}$ in \autoref{def:Dminus}. Clearly,
  \begin{equation}
    \norm{g}_{P_{-}} := \sup_{f\in P_{+}\setminus \set{0}}\frac{\abs{\scprod{g,f}_{\hsmid}}}{\norm{f}_{\hs}}
    \geq \norm{g}_{\hsdual}
  \end{equation}
  and consequently $P_{-} \subseteq D_{-}$.
  Furthermore, we can define $\iota_{P_{+}}$ analogously to $\iota_{+}$. Then we have $\dom \iota_{P_{+}}\adjun = P_{-}$ and $\iota_{+} \subseteq \iota_{P_{+}}$ and therefore $\iota_{P_{+}}\adjun \subseteq \iota_{+}\adjun$.
  Let $f \in P_{+}$. Then there exists a sequence $(f_n)_{n\in\N}$ in $\tilde{D}_{+}$ that converges to $f$ w.r.t. $\norm{.}_{\hs}$. For $g\in P_{-}$ we have
  \begin{align*}
    \abs{\scprod{g,f}_{\hsmid}} &= \abs{\scprod{\iota_{P_{+}}\adjun g, f}_{\hs}}
    = \lim_{n\in\N} \abs{\scprod{\iota_{+}\adjun g, f_{n}}_{\hs}}
    \leq \lim_{n\in\N}\norm{\iota_{+}\adjun g}_{\hs} \norm{f_n}_{\hs}
    = \norm{g}_{\hsdual} \norm{f}_{\hs}
    ,
  \end{align*}
  which yields $\norm{g}_{P_{-}} \leq \norm{g}_{\hsdual}$. Hence, $\norm{.}_{P_{-}} = \norm{.}_{\hsdual}$, $P_{-} = D_{-}$, $\iota_{P_{+}}\adjun = \iota_{+}\adjun$ and $\iota_{P_{+}} = \cl{\iota_{+}}$, which is equivalent to $\hs\cap\hsmid = \cl[\hs\cap\hsmid]{\tilde{D}_{+}}$.
\end{Proof}

We define $D_{+} := \cl[\hs\cap\hsmid]{\tilde{D}_{+}} = \hs\cap\hsmid$ and we will denote the extension of $\iota_{+}$ to $D_{+}$, which is its closure, also by $\iota_{+}$. The adjoint $\iota_{+}\adjun$ is not affected by that.

\begin{Theorem}\label{th:hsdual-hilbert-space}
  The mapping $\iota_{+}\adjun$ can be uniquely extended to a bijective linear isometry $\Psi : \hsdual \to \hs$.
  The space $\hsdual$ is a Hilbert space with the inner product $\scprod{g,f}_{\hsdual} := \scprod{\Psi g, \Psi f}_{\hs}$. Moreover, the induced norm of this inner product coincides with $\norm{.}_{\hsdual}$.
\end{Theorem}

\begin{Proof}
  By \myref{le:adjoint-embedding}{Lemma} $\iota_{+}\adjun$ is a bounded linear mapping from $D_{-}$ to $\hs$ with $\ran \iota_{+}\adjun$ is dense. Since $D_{-}$ is dense in $\hsdual$ by construction, we can extend $\iota_{+}\adjun$ by continuity to $\hsdual$. We denote this extension by $\Psi$. For an arbitrary $g\in \hsdual$ there exists sequence $(g_{n})_{n\in\N}$ in $D_{-}$ that converges to $g$. Hence, 
  \begin{align}\label{eq:Psi-isometry}
    \norm{\Psi g}_{\hs} = \lim_{n\in\N}\norm{\Psi g_{n}}_{\hs} = \lim_{n\in\N} \norm{\iota_{+}\adjun g_{n}}_{\hs}
    = \lim_{n\in\N} \norm{g_{n}}_{\hsdual} = \norm{g}_{\hsdual}.
  \end{align}
  This yields that $\ran \Psi$ is closed in $\hs$. Since $\ran \Psi$ also contains the dense subspace $\ran \iota_{+}\adjun$, the mapping $\Psi$ is surjective.

  Since $\Psi$ is bijective it is easy to see that $\hsdual$ is a Hilbert space with the given inner product. By \eqref{eq:Psi-isometry}, we have
  \begin{equation*}
    \scprod{g,g}_{\hsdual} = \scprod{\Psi g, \Psi g}_{\hs} = \norm{\Psi g}_{\hs}^{2} = \norm{g}_{\hsdual}^{2}.
    \qedhere
  \end{equation*}
\end{Proof}

\begin{Corollary}
  The Hilbert space $\hsdual$ can be identified with the (anti)dual space of $\hs$ by
  \begin{align*}
    \Lambda:\mapping{\hsdual}{\hs\dual}{g}{\scprod{\Psi g, .}_{\hs},}
  \end{align*}
  where $\Psi$ is the mapping from \myref{th:hsdual-hilbert-space}{Theorem}.
\end{Corollary}

\begin{Definition}
  For $f\in\hs$ and $g\in\hsdual$ we define
  \begin{align*}
    \dualprod{g,f}_{\hsdual,\hs} := \dualprod{\Lambda g, f}_{\hs\dual,\hs} = \scprod{\Psi g,f}_{\hs}.
  \end{align*}
  We call $(\hs,\hsmid,\hsdual)$ with this duality a \emph{quasi Gelfand triple}. The space $\hsmid$ will be refered as the \emph{pivot space} and $\Psi$ as the \emph{duality map} in this setting.
\end{Definition}

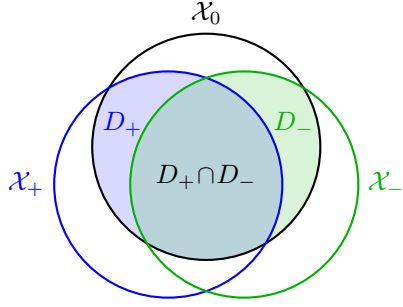
\begin{figure}
\begin{center}
  \begin{tikzpicture}
    \begin{scope}
      \clip (-0.5,-0.5) circle (1.5);
      \fill[blue,opacity=0.15] (0,0) circle (1.5);
    \end{scope}
    \begin{scope}
      \clip (0.5,-0.5) circle (1.5);
      \fill[black!30!green,opacity=0.15] (0,0) circle (1.5);
    \end{scope}
  \draw[thick](0,0) circle [radius=1.5];
  \draw[thick,blue](-0.5,-0.5) circle [radius=1.5];
  \draw[thick,black!30!green](0.5,-0.5) circle [radius=1.5];
  \node[left,blue] at (-2,-0.5) {$\hs$};
  \node[right,black!30!green] at (2,-0.5) {$\hsdual$};
  \node[above] at (0,1.5) {$\hsmid$};
  \node[blue] at (-1.1,0.3) {${D}_{+}$};
  \node[black!30!green] at (1.15,0.3) {$D_{-}$};
    \node at (0,-0.4) {${D}_{+}\!\cap\! D_{-}$};
\end{tikzpicture}
\end{center}
\caption{Illustration of a quasi Gelfand triple}
\label{fig:quasi-gelfand-triple}
\end{figure}

\autoref{fig:quasi-gelfand-triple} illustrates the setting of a quasi Gelfand triple.

\begin{Remark}
  For $f \in D_{+}$ and $g\in D_{-}$ we have
  \begin{align*}
    \dualprod{g,f}_{\hsdual,\hs} = \scprod{\Psi g, f}_{\hs} = \scprod{\iota_{+}\adjun g, f}_{\hs} = \scprod{g,f}_{\hsmid}. 
  \end{align*}
  Since these two sets are dense in $\hs$ and $\hsdual$ respectively, we have for $f\in \hs$ and $g\in\hsdual$
  \begin{align*}
    \dualprod{g,f}_{\hsdual,\hs} = \lim_{(n,m)\in\N^{2}} \scprod{g_{n}, f_{m}}_{\hsmid},
  \end{align*}
  where $(f_{m})_{m\in\N}$ is a sequence in $D_{+}$ that converges to $f$ and $(g_{n})_{n\in\N}$ is a sequence in $D_{-}$ that converges to $g$.
\end{Remark}

In contrast to ``ordinary'' Gelfand triple, the setting for quasi Gelfand triple is somehow ``symmetric'', i.e. the roles of $\hs$ and $\hsdual$ are interchangeable.
If we start with $D_{-}$, then
\begin{equation*}
\iota_{-}: \mapping{D_{-} \subseteq \hsdual}{\hsmid}{g}{g,}
\end{equation*}
is closed by \myref{le:Dminus-complete}{Lemma} and $(D_{-})_{-} = D_{+}$. It is also easy to verify that the unitary operator from $\hs$ to $\hsdual$ resulting from the extension of $\iota_{-}\adjun$ is $\Psi\adjun$. In order to restore $\iota_{-}\adjun$ from $\Psi\adjun$ we only have to restrict $\Psi\adjun$ to $D_{+}$ or more exactly
\begin{equation}\label{eq:iota-minus-adjun}
\iota_{-}\adjun = \Psi\adjun \iota_{+}^{-1}.
\end{equation}

\begin{Proposition}
  The space $D_{+}\cap D_{-}$ is complete with respect to
  $\norm{.}_{\hs\cap\hsdual} := \sqrt{\norm{.}_{\hs}^{2} + \norm{.}_{\hsdual}^{2}}$.
\end{Proposition}

\begin{Proof}
  For $f \in D_{+}\cap D_{-}$ we have
  \begin{align*}
    \norm{f}_{\hsmid}^{2} = \abs{\scprod{f,f}_{\hsmid}} = \abs{\dualprod{f,f}_{\hsdual,\hs}} \leq \norm{f}_{\hsdual} \norm{f}_{\hs} \leq \norm{f}_{\hs\cap\hsdual}^{2}. 
  \end{align*}
  Hence, every Cauchy sequence in $D_{+}\cap D_{-}$ with respect to $\norm{.}_{\hs\cap\hsdual}$ is also a Cauchy sequence with respect to $\norm{.}_{\hsmid}$, $\norm{.}_{\hs}$ and $\norm{.}_{\hsdual}$.

  Let $(f_{n})_{n\in\N}$ be a Cauchy sequencen in $D_{+}\cap D_{-}$ with respect to $\norm{.}_{\hs\cap\hsdual}$. Then the limit with respect to $\norm{.}_{\hsmid}$ and the limit with respect to $\norm{.}_{\hsdual}$ coincide by \myref{le:Dminus-complete}{Lemma}. Furthermore, by the closedness of $\iota_{+}$ the limit with respect to $\norm{.}_{\hsmid}$ and the limit with respect to $\norm{.}_{\hs}$ also coincide. Therefore, all these limits have to coincide and $(f_{n})_{n\in\N}$ converges to that limit in $\norm{.}_{\hs\cap\hsdual}$. 
\end{Proof}

\begin{Lemma}\label{th:iota-p-iota-m-closed}
  The operator
  \begin{align*}
    \begin{bmatrix}
      \iota_{+}    & \iota_{-}
    \end{bmatrix}:
        \mapping{D_{+} \times D_{-} \subseteq \hs \times \hsdual}{\hsmid}{\begin{bmatrix}
          f \\ g
        \end{bmatrix}}{f+g,}    
  \end{align*}
  is closed.
\end{Lemma}

\begin{Proof}
  Let $\left(\left(\begin{bsmallmatrix}f_n \\ g_n\end{bsmallmatrix},z_n\right)\right)_{n\in\N}$ be a sequence in $\begin{bmatrix}\iota_{+}    & \iota_{-} \end{bmatrix}$ that converges to $\left(\begin{bsmallmatrix}f \\ g\end{bsmallmatrix},z\right) \in \hs\times \hsdual \times \hsmid$. Then we have
    \begin{align*}
      \norm{z}_{\hsmid}^2 = \lim_{n\in\N} \norm{f_n + g_n}_{\hsmid}^2 = \lim_{n\in\N}\big( \norm{f_n}_{\hsmid}^2 + \norm{g_n}_{\hsmid}^2 + 2\Re \scprod{f_n,g_n}_{\hsmid}\big).
    \end{align*}
    Since $2\Re \scprod{f_n,g_n}_{\hsmid}$ converges to $2\Re\dualprod{f,g}_{\hs,\hsdual}$, we conclude that $\norm{f_n}_{\hsmid}$ and $\norm{g_n}_{\hsmid}$ are bounded. Hence, it exists a subsequence $(f_{n(k)})_{k\in\N}$ that converges weakly to an $\tilde{f}\in\hsmid$. Moreover, by \autoref{th:weak-to-strong-convergent} there is a further subsequence such that $\frac{1}{j}\sum_{i=1}^{j}f_{n(k(i))}$ converges to $\tilde{f}$ strongly. The sequence $\big(\frac{1}{j}\sum_{i=1}^{j}f_{n(k(i))}\big)_{j\in\N}$ has still the limit $f$ in $\hs$ and because $\iota_{+}$ is closed we conclude that $f=\tilde{f}\in D_{+}$. We also have $\frac{1}{j}\sum_{i=1}^{j}g_{n(k(j))} \to z - f$ in $\hsmid$. By \myref{le:Dminus-complete}{Lemma} $g = z-f \in D_{-}$. Hence, the operator $\begin{bmatrix}\iota_{+} & \iota_{-}\end{bmatrix}$ is closed.
\end{Proof}

\begin{Proposition}
  $D_{+} \cap D_{-}$ is dense in $\hsmid$ with respect to $\norm{.}_{\hsmid}$.
\end{Proposition}

\begin{Proof}
  By $\dom \iota_{\pm}\adjun = D_{\mp}$ we have
  \begin{equation*}
    \hsmid = \big(\mul \begin{bmatrix} \iota_{+} & \iota_{-}\end{bmatrix}\big)^{\perp} = \cl{\dom \begin{bmatrix} \iota_{+} & \iota_{-}\end{bmatrix}\adjun} = \cl{\dom \iota_{+}\adjun \cap \dom \iota_{-}\adjun} = \cl{D_{-} \cap D_{+}}
    .
    \qedhere
  \end{equation*}
\end{Proof}

The following theorem can be found in \cite[Theorem 2 p. 200]{fana-yosida}, we just changed that the operator maps into a different space. Hence, we provide a proof.

\begin{Theorem}[J. von Neumann]\label{th:TTadjun-self-adjoint}
Let $T$ be a closed linear operator from the Hilbert spaces $X$ to the Hilbert space $Y$. Then $T\adjun T$ and $TT\adjun$ are self-adjoint, and $(\opid_{X} + T\adjun T)$ and $(\opid_{Y} + TT\adjun)$ are boundedly invertible.
\end{Theorem}

\begin{Proof}
Since $T\adjun = \left[\begin{smallmatrix} 0 & \opid_{Y} \\ -\opid_{X} & 0\end{smallmatrix}\right] T^{\perp}$, we have $T \oplus \left[\begin{smallmatrix} 0 & -\opid_{X} \\ \opid_{Y} & 0\end{smallmatrix}\right] T\adjun = X\times Y$. Hence, for $(h,0)\in X\times Y$ there are unique $x \in \dom T$ and $y \in \dom T\adjun$ such that
\begin{equation}\label{eq:decomposition}
(h,0) = (x,Tx) + (-T\adjun y,y).
\end{equation}
Consequently, $h=x - T\adjun y$ and $y = -Tx$, which implies $x \in \dom T\adjun T$ and
\begin{align*}
h = x + T\adjun T x.
\end{align*}
Because of the uniqueness of the decomposition in \eqref{eq:decomposition}, $x \in \dom T\adjun T$ is uniquely determined by $h \in X$. Therefore, $(\opid_{X} + T\adjun T)^{-1}$ is a well-defined and everywhere defined operator.

For $h_{1},h_{2} \in X$, we define $x_{1} := (\opid_{X} + T\adjun T)^{-1} h_{1}$ and $x_{2} := (\opid_{X} + T\adjun T)^{-1} h_{2}$. Then $x_{1}, x_{2} \in \dom T\adjun T$ and, by the closedness  of $T$, $T\adjun[2] = T$. Hence,
\begin{align*}
\scprod{h_{1}, (\opid_{X} + T\adjun T)^{-1} h_{2}} &= \scprod{(\opid_{X} + T\adjun T) x_{1}, x_{2}} = \scprod{x_{1},x_{2}} + \scprod{T\adjun T x_{1}, x_{2}} \\
&= \scprod{x_{1},x_{2}} + \scprod{T x_{1},T x_{2}}= \scprod{x_{1},x_{2}} + \scprod{x_{1}, T\adjun T x_{2}} \\
&= \scprod{x_{1}, (\opid_{X} + T\adjun T) x_{2}} = \scprod{(\opid_{X} + T\adjun T)^{-1} h_{1}, h_{2}},
\end{align*}
which yields that $(\opid_{X} + T\adjun T)^{-1}$ is self-adjoint. Therefore $(\opid_{X} + T\adjun T)$ and $T\adjun T$ are also self-adjoint. Moreover, $(\opid_{X} + T\adjun T)^{-1}$ is bounded as a closed and everywhere defined operator.

By $TT\adjun = (T\adjun)\adjun(T\adjun)$ the other statements follow by the already shown.
\end{Proof}

Applying this theorem to $S=\lambda T$ implies that $\R_{-}$ is contained in the resolvent set of $T\adjun T$.

\begin{Corollary}
  The set $D_{+}\cap D_{-}$ is dense in $\hs$ and $\hsdual$ with respect to their corresponding norms.
\end{Corollary}

\begin{Proof}
  Applying \myref{th:TTadjun-self-adjoint}{Theorem} to $\iota_{+}$ yields $\iota_{+}\adjun \iota_{+}$ is self-adjoint. Hence, $\dom \iota_{+}\adjun \iota_{+}$ is dense in $\hs$. By \myref{le:adjoint-embedding}{Lemma} $\dom \iota_{+}\adjun = D_{-}$, consequently $\dom \iota_{+}\adjun \iota_{+} = D_{+}\cap D_{-}$.

  An analogous argument for $\iota_{-}$ yields $D_{+}\cap D_{-}$ is dense in $\hsdual$.
\end{Proof}

\begin{Corollary}
  $D_{+} + D_{-} = \hsmid$.
\end{Corollary}

\begin{Proof}
  Applying \myref{th:TTadjun-self-adjoint}{Theorem} to $\iota_{+}$ yields $(\opid_{\hsmid} + \iota_{+}\iota_{+}\adjun)$ is onto. Hence, for every $x\in \hsmid$ there exists a $g_{x} \in \dom \iota_{+}\iota_{+}\adjun \subseteq D_{-}$ such that
  \begin{equation*}
    x = \underbrace{g_x}_{\in D_{-}} + \underbrace{\iota_{+} \iota_{+}\adjun g_x}_{\in D_{+}}
    .
  \end{equation*}
  Since  $g_{x} \in \dom \iota_{+}\iota_{+}\adjun$, we have $\iota_{+}\adjun g_x \in D_{+}$ and consequently $x \in D_{+} + D_{-}$.
\end{Proof}

\begin{Proposition}\label{prop:TDplus-new-triple}
Let $T$ be a bounded and boundedly invertible mapping on $\hsmid$. Then $P_{+} := T D_{+}$ equipped with $\norm{f}_{\mathcal{Y}_{+}} := \norm{T^{-1}f}_{\hs}$ establishes a quasi Gelfand triple $(\mathcal{Y}_{+},\hsmid, \mathcal{Y}_{-})$, where $\mathcal{Y}_{+}$ is the completion of $P_{+}$ and $\mathcal{Y}_{-}$ is the completion of $P_{-}$ defined as in \myref{def:Dminus}{Definition} where $D_{+}$ is replaced by $P_{+}$. Moreover, $P_{-} = (T\adjun)^{-1} D_{-}$, $\norm{g}_{\mathcal{Y}_{-}} = \norm{T\adjun g}_{\hsdual}$, and $T$ and $(T\adjun)^{-1}$ can be continuously extended to linear bounded and boundedly invertible mappings from $\hs$ and $\hsdual$ to $\mc Y_{+}$ and $\mc Y_{-}$ respectively.
\end{Proposition}

\begin{Proof}
The mapping $T\big\vert_{D_{+}}$ is also bounded and boundedly invertible if we equipped its domain with $\norm{.}_{\hs}$ and its codomain with $\norm{.}_{\mathcal{Y}_{+}}$. 
So the linear relation $\big[\begin{smallmatrix} T & 0 \\ 0 & T\end{smallmatrix}\big] \iota_{+} = \set{(Tf,Tg) : (f,g) \in \iota_{+}} \subseteq \mc{Y}_{+} \times \hsmid$ is closed. Since this linear relation coincides with the embedding $\iota_{P_{+}}: P_{+} \subseteq \mathcal{Y}_{+} \to \hsmid, f \mapsto f$, \myref{prop:iota-closed-equivalences}{Proposition} yields that all assumptions for a quasi Gelfand triple are satisfied. For $g \in \hsmid$ we have
\begin{equation*}
\norm{g}_{\mathcal{Y}_{-}} = \sup_{h\in P_{+}} \frac{\abs{\scprod{g,h}_{\hsmid}}}{\norm{h}_{\mathcal{Y}_{+}}} = \sup_{f\in D_{+}} \frac{\abs{\scprod{g,Tf}_{\hsmid}}}{\norm{Tf}_{\mathcal{Y}_{+}}}
= \sup_{f\in D_{+}} \frac{\abs{\scprod{T\adjun g,f}_{\hsmid}}}{\norm{f}_{\hs}} = \norm{T\adjun g}_{\hsdual}.\qedhere
\end{equation*}
\end{Proof}

\begin{Corollary}\label{co:closed-operator-from-Dp-times-Dm}
Let $S,T$ be a bounded and boundedly invertible mappings on $\hsmid$. Then $\big[\begin{smallmatrix} ST\big\vert_{D_{+}} & S(T\adjun)^{-1}\big\vert_{D_{-}}\end{smallmatrix}\big]$ is a closed linear relation between $\hs\times \hsdual$ and $\hsmid$ with 
$\ran \big[\begin{smallmatrix} ST\big\vert_{D_{+}} & S(T\adjun)^{-1}\big\vert_{D_{-}}\end{smallmatrix}\big] = \hsmid$.
\end{Corollary}

\begin{Proof}
  Let $P_{+} = T D_{+}$. Then by \myref{prop:TDplus-new-triple}{Proposition} the corresponding $P_{-}$ can be obtained by $(T\adjun)^{-1} D_{-}$. The mapping
  \begin{equation*}
    \Xi: \mapping{\hs \times \hsdual \times\hsmid}
    {\mc Y_{+} \times \mc Y_{-} \times \hsmid}
    {
      \left[
      \begin{matrix}
        f \\ g \\ z
      \end{matrix}
    \right]
    }{
    \left[
      \begin{matrix}
        T & 0 & 0 \\
        0 & (T\adjun)^{-1} & 0 \\
        0 & 0 & S^{-1}
      \end{matrix}
    \right]
    \left[
      \begin{matrix}
        f \\ g \\ z
      \end{matrix}
    \right]
  }
  \end{equation*}
  is linear bounded and boundedly invertible, where $\mc Y_{+}$ and $\mc{Y}_{-}$ are the spaces corresponding to $P_{+}$ and $P_{-}$ from \myref{prop:TDplus-new-triple}{Proposition}.
  Since $(\mc Y_+, \hsmid, \mc Y_-)$ is a quasi Gelfand triple,
  $\left[
    \begin{matrix}
      \iota_{P_{+}} & \iota_{P_{-}}
    \end{matrix}
  \right]$
  is closed in $\mc Y_+ \times \mc Y_- \times \hsmid$ and therefore also its pre-image under $\Phi$
  \begin{equation*}
    \Xi^{-1}\left(
      \left[
        \begin{matrix}
          \iota_{P_{+}} & \iota_{P_{-}}
        \end{matrix}
      \right]
    \right)
    =
    \left[
      \begin{matrix}
        T^{-1} & 0 & 0 \\
        0 & T\adjun & 0 \\
        0 & 0 & S
      \end{matrix}
    \right]
    \left[
      \begin{matrix}
        \iota_{P_{+}} & \iota_{P_{-}}
      \end{matrix}
    \right]
    =
    \left[
      \begin{matrix}
        ST\iota_{{+}} & S(T\adjun)^{-1}\iota_{{-}}
      \end{matrix}
    \right]
  \end{equation*}
  is closed. Furthermore,
  \begin{equation*}
    \ran \left[
      \begin{matrix}
        ST\big\vert_{D_{+}} & S(T\adjun)^{-1}\big\vert_{D_{-}}
      \end{matrix}
    \right]
    = S \ran
    \left[
      \begin{matrix}
        \iota_{P_{+}} & \iota_{P_{-}}
      \end{matrix}
    \right]
    =
    S \hsmid = \hsmid
    .
    \qedhere
  \end{equation*}
\end{Proof}

\begin{Lemma}\label{th:trans-boundary-triple}
  Let $T$ be a bounded and boundedly invertible mapping on $\hsmid$ and $(\hs, \hsmid, \hsdual)$ be a quasi Gelfand triple such that $(\hs, B_{1}, \Psi B_{2})$ is a boundary triple for an operator $A$. Furthermore, let $\mc Y_{+}$ and $\mc Y_{-}$ be as defined in \autoref{prop:TDplus-new-triple}.
  Then $(\mc Y_{+}, \hsmid, \mc Y_{-})$ is also a quasi Gelfand triple such that $(\mc Y_{+}, TB_{1}, \Phi (T\adjun)^{-1}B_{2})$ is a boundary triple for $A$, where $\Phi$ denotes the duality map of $(\mc Y_{+}, \hsmid, \mc Y_{-})$.
\end{Lemma}

\begin{Proof}
  By \autoref{prop:TDplus-new-triple} $(\mc Y_{+}, \hsmid, \mc Y_{-})$ is a quasi Gelfand triple. Note that $T$ and $(T\adjun)^{-1}$ can be extended to mappings from $\hs$ to $\mc Y_{+}$ and $\hsdual$ to $\mc Y_{-}$ respectively. For $x,y \in \dom A$ we have
  \begin{align*}
    \scprod{B_{1}x, \Psi B_{2} y}_{\hs} &= \dualprod{B_{1}x,B_{2}y}_{\hs,\hsdual} = \dualprod{TB_{1}x, (T\adjun)^{-1}B_{2}y}_{\mc Y_{+}, \mc Y_{-}}
    \\ &=  \scprod{TB_{1}x, \Phi (T\adjun)^{-1}B_{2} y}_{\mc Y_{+}}
    .
    \qedhere
  \end{align*}
\end{Proof}

\section{Boundary Spaces}

In this section we will construct a suitable boundary space $\PLspace$ (\autoref{def:boundary-space}), where we will later formulate boundary conditions.
This space will provide a quasi Gelfand triple with a subspace of $\Lp{2}(\partial \Omega)$ as pivot space.
\smallskip

\begin{Definition}\label{def:splitting-with-thin-boundaries}
  We say $(\Gamma_{j})_{j=1}^{k}$, where $\Gamma_{j} \subseteq \partial \Omega$, is a \emph{splitting with thin boundaries} of $\partial \Omega$, if
  \begin{enumerate}[label = \textrm{\textup{(\roman*)}}]
    \item $\bigcup_{j=1}^{k} \cl{\Gamma_{j}} = \partial \Omega$,
    
    \item the sets $\Gamma_{j}$ are pairwise disjoint,
    
    \item the sets $\Gamma_{j}$ are relatively open in $\partial \Omega$,
    
    \item the boundaries of $\Gamma_{j}$ have zero measure w.r.t. the surface measure of $\partial \Omega$.
  
  \end{enumerate}
\end{Definition}

For $\Gamma \subseteq \partial \Omega$ we will denote by $P_{\Gamma}$ the orthogonal projection from $\Lp{2}(\partial\Omega)^{\irow}$ on $\Lp{2}_{\pi}(\Gamma) := \cl{\ran \indicator_{\Gamma}\Lnub} \leq \Lp{2}(\Gamma)^{\irow}$. Therefore, we can adapt \eqref{eq:int-by-parts} such that
\begin{align}\label{eq:int-by-parts-projection}
	 \scprod{\diffop f,g}_{\Lp{2}(\Omega)^{\irow}} + \scprod{f,\diffopad g}_{\Lp{2}(\Omega)^{\icol}}
    	= \scprod{\Lnub \boundtr f,\underbrace{P_{\partial\Omega}\boundtr}_{\displaystyle\pi_{\myL}} g}_{\Lp{2}(\partial \Omega)^{\irow}}
	.
\end{align}
We define $\pi_{\myL}^{\Gamma} : H^{1}(\Omega)^{\irow} \to \Lp{2}_{\pi}(\Gamma)\leq \Lp{2}(\Gamma)^{\irow}$ by $\pi_{\myL}^{\Gamma} := P_{\Gamma}\boundtr$ and $\pi_{\myL} := \pi_{\myL}^{\partial \Omega}$. Since both $P_{\Gamma}$ and $\boundtr$ are continuous, the mapping $\pi_{\myL}^{\Gamma}$ is also continuous. Therefore, $\ker \pi_{\myL}^{\Gamma}$ is closed. Note that $P_{\Gamma} = \indicator_{\Gamma} P_{\partial \Omega}$ and consequently $\pi_{\myL}^{\Gamma} = \indicator_{\Gamma} \pi_{\myL}$, and $\indicator_{\Gamma} \Lnub = \Lnub \indicator_{\Gamma}$.

\begin{Example}\label{ex:div-grad-2}
  Let $\myL$ be as in \myref{ex:div-grad}{Example}. Then $\Lnub f = \scprod{f,\nu}_{\K^{3}}$ and $\Lnub$ is certainly surjective. Therefore, $\Lp{2}_{\pi}(\partial \Omega) = \Lp{2}(\partial \Omega)$, $\piL = \boundtr$ and $\piL[\Gamma] = \indicator_{\Gamma}\boundtr$. Since $\diffopad = \grad$, we have that $\Ladspace = H^{1}\mspace{-1.5mu}(\Omega)$.
\end{Example}

\begin{Lemma}\label{le:ker-closed}
  Let $\Gamma \subseteq \partial \Omega$ be relatively open. Then the subspace $\ker  \pi_{\myL}^{\Gamma}$ is closed in $H^{1}(\Omega)^{\irow}$ with respect to $\Ladnorm{.}$. This can also be formulated as
  \begin{align*}
    \cl[\Ladnorm{.}]{\ker \pi_{\myL}^{\Gamma}} \cap H^{1}\mspace{-1.5mu}(\Omega)^{\irow} = \ker \pi_{\myL}^{\Gamma}
    .
  \end{align*}
\end{Lemma}

\begin{Proof}
  Let $(g_{n})_{n\in\N}$ be a sequence in $\ker \pi_{\myL}^{\Gamma}$ which converges to $g \in H^{1}(\Omega)^{\irow}$ with respect to $\Ladnorm{.}$. By \myref{co:int-by-parts}{Corollary} we have for an arbitrary $f \in H^{1}_{\partial \Omega \setminus \Gamma}(\Omega)^{\icol} := \set{f \in H^{1}(\Omega)^{\icol} : \indicator_{\partial \Omega \setminus \Gamma}\boundtr f = 0}$
  \begin{align*}
    \abs{\scprod{\Lnub \boundtr f, \pi_{\myL}^{\Gamma} (g-g_{n})}_{\Lp{2}} }
    = \abs{\scprod{\Lnub \boundtr f, \pi_{\myL} (g-g_{n})}_{\Lp{2}} }
    \leq \Lnorm{f} \Ladnorm{g-g_{n}}
    .
  \end{align*}
  Since $\pi_{\myL}^{\Gamma} (g-g_{n}) = \pi_{\myL}^{\Gamma} g$ and the right-hand-side converges to $0$, we can see that $\pi_{\myL}^{\Gamma} g \perp \Lnub \boundtr H^{1}_{\partial\Omega\setminus\Gamma}(\Omega)^{\icol}$.
  By \cite[Th. 13.6.10, Re. 13.6.12]{observation-and-control} $\boundtr H^{1}_{\partial \Omega \setminus \Gamma}(\Omega)^{\icol}$ is dense in $\Lp{2}(\Gamma)^{\icol}$, which implies $\pi_{\myL}^{\Gamma} g \perp \ran \indicator_{\Gamma} \Lnub$. By definition $\pi_{\myL}^{\Gamma} g$ is also in $\cl{\ran \indicator_{\Gamma}\Lnub}$, which leads to $\pi_{\myL}^{\Gamma} g = 0$.
  Hence, $\ker \pi_{\myL}^{\Gamma}$ is closed in $H^{1}(\Omega)^{\irow}$ with respect to $\Ladnorm{.}$.
\end{Proof}

By the previous lemma we can endow $M_{\Gamma} := \ran \pi_{\myL}^{\Gamma}$ with the norm
\begin{align*}
  \norm{\phi}_{M_{\Gamma}} := \inf\set*{\Ladnorm{g} : \pi_{\myL}^{\Gamma} g = \phi}
  ,
\end{align*}
which makes it a pre-Hilbert space. The next lemma will clarify that.

\begin{Lemma}\label{le:factor-spaces}
  The space $(M_{\Gamma}, \norm{.}_{M_{\Gamma}})$ is a pre-Hilbert space. Furthermore, its completion denoted by $(\cl{M_{\Gamma}},\norm{.}_{\cl{M_{\Gamma}}})$ is isomorphic to $\quspace{\Ladspace}{\cl[\Ladspace]{\ker \piL[\Gamma]}}$. The mapping $\piL[\Gamma]$ can be continuously extended to $\Ladspace$.
  For the kernel of the extension $\cl{\piL[\Gamma]}$ we have $\ker \cl{\piL[\Gamma]} = \cl[\Ladspace]{\ker \piL[\Gamma]}$.
\end{Lemma}

\begin{Proof}
  \renewcommand{\piL}{\pi_{\myL}^{\Gamma}}%
  By \myref{le:ker-closed}{Lemma} $\ker \piL$ is closed in $H^{1}(\Omega)^{\irow}$ with respect to $\Ladnorm{.}$, which implies that $\Big( \quspace{H^{1}(\Omega)^{\irow}}{\ker \piL},\norm{.}_{ \quspace[\scriptstyle]{\Ladspace\;}{\ker \piL}} \Big)$ is a normed space. Since
  \begin{align*}
    \norm*{[g]_{\raisebox{-2pt}{$\scriptstyle\sim$}}}_{ \quspace[\scriptstyle]{\Ladspace\;}{\ker \piL}}
    =
    \norm*{\piL g}_{M_{\Gamma}}
    ,
  \end{align*}
  it is straight forward that $(M_{\Gamma},\norm{.}_{M_{\Gamma}})$ is isomorphic to $\Big( \quspace{H^{1}(\Omega)^{\irow}}{\ker \piL},\linebreak[4]\norm{.}_{ \quspace[\scriptstyle]{\Ladspace\;}{\ker \piL}} \Big)$.
  
  Clearly, $(M_{\Gamma},\norm{.}_{M_{\Gamma}})$ has a completion $(\cl{M_{\Gamma}},\norm{.}_{\cl{M_{\Gamma}}})$. By Definition of the norm $\norm{.}_{M_{\Gamma}}$ we have for every $g \in H^{1}(\Omega)^{\irow}$
  \begin{align*}
    \norm{\piL g}_{\cl{M_{\Gamma}}} = \norm{\piL g}_{M_{\Gamma}} \leq \Ladnorm{g}
    .
  \end{align*}
  Therefore, we can extend $\piL$ by continuity on $\Ladspace$. To avoid confusion in this proof we will use the symbol  $\cl{\piL}$ for this extension.

  Let $g\in \Ladspace$ and $(g_{n})_{n\in\N}$ a sequence in $H^{1}(\Omega)^{\irow}$ which converges to $g$. Then we have
  \begin{align*}
    \norm{\cl{\piL} g}_{\cl{M_{\Gamma}}} = \lim_{n\in\N}\norm{\piL g_{n}}_{{M_{\Gamma}}}
    = \lim_{n\in\N}\inf_{k\in\ker \piL} \Ladnorm{g_{n} + k}
  \end{align*}
  The triangular inequality yields
  \begin{align*}
    \inf_{k\in \ker \piL}\norm{g + k} - \norm{g_{n} - g}
    \leq \inf_{k\in\ker \piL} \norm{g_{n} + k}
    \leq \inf_{k\in \ker \piL}\norm{g + k} + \norm{g_{n} - g}
    .
  \end{align*}
  Hence, we have $\norm{\cl{\piL} g}_{\cl{M_{\Gamma}}} = \inf_{k\in \ker \piL}\norm{g + k} = \inf_{k\in \cl{\ker \piL}}\norm{g + k}$ and consequently $\quspace{\Ladspace}{\cl{\ker \piL}}$ is isomorphic to $\ran \cl{\piL}$. Since $\quspace{\Ladspace}{\cl{\ker \piL}}$ is a Hilbert space, in particular complete, and $M_{\Gamma} \subseteq \ran \cl{\piL} \subseteq \cl{M_{\Gamma}}$, we have $\cl{M_{\Gamma}} = \ran \cl{\piL}$, which makes $\cl{M_{\Gamma}}$ also a Hilbert space.
\end{Proof}

We will use the symbol $\piL[\Gamma]$ also for its continuous extension $\cl{\piL[\Gamma]}$.

\begin{Definition}\label{def:boundary-space}
  Let $\Gamma_0, \Gamma_1 \subseteq \partial \Omega$ be a splitting with thin boundaries. Then we define
\begin{equation*}
  \Ladspace[\Gamma_0] := \ker \pi_{\myL}^{\Gamma_0}
  \quad \text{and} \quad
  \PLspace[\Gamma_1] := \ran \pi_{\myL}\big\vert_{\Ladspace[\Gamma_0]}
  ,
\end{equation*}
where we endow $\Ladspace[\Gamma_0]$ with $\Ladnorm{.}$ and $\PLspace[\Gamma_1]$ with $\norm{.}_{\PLspace[\Gamma_1]} := \norm{.}_{\cl{M_{\partial\Omega}}}$. Instead of $\PLspace[\partial \Omega] = \ran \piL = \cl{M_{\partial \Omega}}$ we just write $\PLspace$.
\end{Definition}

\begin{Example}\label{ex:div-grad-3}
  Continuing \autoref{ex:div-grad-2} yields
  $\Ladspace[\Gamma_{0}] = H^{1}_{\Gamma_{0}}(\Omega)^{\irow} = \set{f \in H^{1}(\Omega)^{\irow} : \indicator_{\Gamma_{1}}\boundtr f = 0}$ which already appeared in the proof of \autoref{le:ker-closed}. Moreover, we have $\piL = \boundtr$, $\piL[\Gamma_{1}] = \indicator_{\Gamma_{1}}\boundtr$,
  $\PLspace = H^{\nicefrac{1}{2}}(\partial\Omega)$, and
  $\PLspace[\Gamma_{1}] = \set{f \in H^{\nicefrac{1}{2}}(\partial \Omega) : f\big\vert_{\Gamma_{0}} = 0}$.
\end{Example}

\begin{Remark}\label{re:boundary-spaces-dense}
  $\Ladspace[\Gamma_0]$ is again a Hilbert space and $H^1(\Omega)^{\irow} \cap \Ladspace[\Gamma_0]$ is dense in $\Ladspace[\Gamma_0]$. The density follows from the assertion $\cl{\ker \piL[\Gamma]\big\vert_{H^1(\Omega)^{\irow}}} = \ker \piL[\Gamma]$ of the previous lemma.
Moreover, $\PLspace[\Gamma_1]$ is closed in $\PLspace$, since $\piL[\Gamma_0]\circ\piL^{-1}$ is well-defined and continuous, and $\PLspace[\Gamma_1] = \ker \piL[\Gamma_0] \circ \piL^{-1}$. Hence, $\PLspace[\Gamma_1]$ is also a Hilbert space.
\end{Remark}

For $g \in H^{1}(\Omega)^{\irow} \cap \Ladspace[\Gamma_{0}]$, we have that $\piL g = \piL[\Gamma_{1}] g$ as elements of $\Lp{2}(\partial \Omega)$. So somehow it is possible to say that $\PLspace[\Gamma_{1}] = \cl{M_{\Gamma_{1}}}$, but the norms are different.

\begin{Proposition}\label{th:Lnu-maps-into-dualspace}
  The mapping $\indicator_{\Gamma_{1}}\mspace{-1mu}\Lnub \boundtr : H^{1}(\Omega)^{\icol} \to \Lp{2}_{\pi}(\Gamma_{1})$ can be extended to a 
  linear continuous mapping
  \begin{align*}
    \Lnu[\Gamma_1] : \Lspace \to \PLspace[\Gamma_{1}]\dual
    ,
  \end{align*} 
  such that $\norm{\Lnu[\Gamma_1] f}_{\PLspace[\Gamma_1]\dual} \leq  \Lnorm{f}$.
\end{Proposition}

\begin{Proof}
  Let $f \in H^{1}(\Omega)^{\icol}$. For $g\in H^{1}(\Omega)^{\irow}\cap \Ladspace[\Gamma_0]$ we have by \myref{co:int-by-parts}{Corollary}
  \begin{align*}
    \abs*{\scprod{\indicator_{\Gamma_1}\Lnub \boundtr f,\piL g}_{\Lp{2}(\Gamma_1)^{\irow}}}
    = \abs*{\scprod{\Lnub \boundtr f, \piL g}_{\Lp{2}(\partial \Omega)^{\irow}}}
    \leq \Lnorm{f} \Ladnorm{g}
    .
  \end{align*}
  By \myref{re:boundary-spaces-dense}{Remark} it is easy to see that the subspace $M := \ran \piL\big\vert_{H^{1}(\Omega)^{\irow}\cap \Ladspace[\Gamma_0]} \subseteq \Lp{2}_{\pi}(\Gamma_1)^{\irow}$ of $\PLspace[\Gamma_1]$ is dense.
  For $\phi \in M$ there exists at least one $g\in H^{1}(\Omega)^{\irow}\cap \Ladspace[\Gamma_0]$ such that $\pi_{\myL} g = \phi$. Hence, we can rewrite the inequality by
  \begin{align*}
    \abs*{\scprod{\indicator_{\Gamma_1}\Lnub \boundtr f, \phi}_{\Lp{2}(\Gamma_1)^{\irow}}}
    &\leq \Lnorm{f} \inf_{\substack{g\in H^{1}(\Omega)^{\irow} \cap \Ladspace[\Gamma_{0}] \\ \pi_{\myL}g = \phi}}\Ladnorm{g} \\
    &= \Lnorm{f} \norm{\phi}_{\PLspace[\Gamma_1]}
      .
  \end{align*}
  We will extend the mapping $\phi \mapsto \scprod{\indicator_{\Gamma_1}\Lnub \boundtr f, \phi}_{\Lp{2}(\Gamma_1)^{\irow}}$ by continuity on $\PLspace[\Gamma_{1}]$. We will denote this extension by $\Xi_{f}$.
  Therefore, we have
  \begin{align*}
    \abs*{
    \Xi_{f}(\phi)
    }
    \leq \Lnorm{f} \norm{\phi}_{\PLspace[\Gamma_1]}
    .
  \end{align*}
  This means that the mapping
  $f \mapsto
  \Xi_{f}
  $
  is continuous, if we endow $H^{1}(\Omega)^{\icol}$ with $\Lnorm{.}$. Once again, we will extend this mapping by continuity on $\Lspace$ and denote it by $\Lnu[\Gamma_1]$.
\end{Proof}

Instead of writing $\Lnu[\partial\Omega\mspace{2mu}]$ we will just write $\Lnu$.

\begin{Remark}
  In fact the extension of the $\Lp{2}(\Gamma_1)$ scalar product in the previous proof is nothing else but 
  \begin{align*}
    (\Lnu[\Gamma_{1}] f)(\phi)
    =
    \scprod{\diffop f,g}_{\Lp{2}(\Omega)^{\irow}} + \scprod{f,\diffopad g}_{\Lp{2}(\Omega)^{\icol}}
    ,
  \end{align*}
  where $g \in \Ladspace[\Gamma_0]$ is any element that satisfies $\pi_{\myL} g = \phi$.
\end{Remark}

\begin{Remark}\label{re:Lnu-equals-LnuGamma}
Since $\PLspace[\Gamma_{1}]$ is a subspace of $\PLspace[\partial \Omega] = \PLspace$ every element of $\PLspace\dual$ can also be treated as an element of $\PLspace[\Gamma_{1}]\dual$. By definition of $\Lnu[\Gamma_{1}]$ and $\Lnu$ it is easy to see that $\Lnu[\Gamma_{1}] f = \Lnu f\big\vert_{\PLspace[\Gamma_{1}]}$ or equivalently $\Lnu[\Gamma_{1}] f$ and $\Lnu f$ coincide as elements of $\PLspace[\Gamma_{1}]\dual$. The reason for even defining $\Lnu[\Gamma_{1}]$ is that the range of its restriction to $H^{1}(\Omega)^{\icol}$ is also contained in $\Lp{2}_{\pi}(\Gamma_{1})$, which will be important for getting a quasi Gelfand triple.
\end{Remark}

\begin{Corollary}\label{co:int-by-parts-extended}
  For $f \in \Lspace$ and $g \in \Ladspace$ we have
  \begin{align*}
    \scprod{\diffop f,g}_{\Lp{2}(\Omega)^{\irow}} + \scprod{f,\diffopad g}_{\Lp{2}(\Omega)^{\icol}}
    &
      = \dualprod{\Lnu f,\piL g}_{\PLspace\dual, \PLspace} 
    \\
    & = \dualprod{\pi_{\myL\hermitian} f, \Lnuad g}_{\PLadspace, \PLadspace\dual}
      ,
  \end{align*}
  and for $h \in \Ladspace$ such that $\diffop \diffopad h \in \Lp{2}(\Omega)^{\irow}$ we have
  \begin{equation}\label{eq:green-id}
    \scprod{\diffop \diffopad h,g}_{\Lp{2}(\Omega)^{\irow}} + \scprod{\diffopad h,\diffopad g}_{\Lp{2}(\Omega)^{\icol}}
    = \dualprod{\Lnu \diffopad h,\piL g}_{\PLspace\dual, \PLspace}
    .
  \end{equation}
\end{Corollary}

\begin{Proof}
  Since $H^{1}(\Omega)^{\icol}$ is dense in $\Lspace$ and
  $H^{1}(\Omega)^{\irow}$ is dense in $\Ladspace$, the first equation follows from \eqref{eq:int-by-parts-projection} by continuity.
  Switching the roles of $\diffop$ and $\diffopad$ yields the second equation.

  For the second assertion set $f=\diffopad h$ in the first equation.
\end{Proof}

\begin{Remark}
For $g \in \Ladspace[0]$ there is a sequence $(g_{n})_{n\in\N}$ in $\mc D(\Omega)$ converging to $g$, which yields $\piL g = \lim_{n\in\N} \piL g_{n} = 0$. Therefore,
	$\Ladspace[0] \subseteq  \ker \piL = \Ladspace[\partial \Omega]$. On the other hand, if $g \in \Ladspace[\partial \Omega]$, then
	\begin{equation*}
		\scprod{\diffop f,g}_{\Lp{2}(\Omega)^{\irow}} + \scprod{f,\diffopad g}_{\Lp{2}(\Omega)^{\icol}}
      		= \dualprod{\Lnu f,\piL g}_{\PLspace\dual, \PLspace} = 0
		.
	\end{equation*}
	Hence, by \myref{le:Lspace0}{Lemma} $g \in \Ladspace[0]$. Consequently, $\Ladspace[0] = \Ladspace[\partial \Omega]$. Clearly the same holds true for $\Lspace$.
\end{Remark}

\begin{Theorem}\label{th:Lnu-onto}
  The mapping $\Lnu : \Lspace \to \PLspace\dual$ is 
  linear, bounded and onto.
\end{Theorem}

\begin{Proof}
  \newcommand{\dualsc}[2][]{\dualprod[#1]{#2}_{\PLspace\dual, \PLspace}}%
  \newcommand{\testdual}[2][]{\dualprod[#1]{#2}_{\mc D\dual(\Omega)^{}, \mc D(\Omega)}}%
  By \myref{th:Lnu-maps-into-dualspace}{Proposition} we already know that $\Lnu$ is 
  linear and bounded, and maps $\Lspace$ into $\PLspace\dual$.
  
  Let $\mu \in \PLspace\dual$ be arbitrary. Since $\pi_{\myL}$ is continuous from $\Ladspace$ to $\PLspace$, the mapping $g \mapsto \dualsc{ \mu, \piL g}$ is also continuous. Consequently, there exists an $h \in \Ladspace$ such that
  \begin{align*}
    \scprod{h,g}_{\Ladspace} = \dualsc{\mu, \pi_{\myL}g}
    \quad
    \text{for all}\quad g\in \Ladspace
    .
  \end{align*}
  For  a test function $v \in \mathcal{D}(\Omega)^{\icol}$ we have
  \begin{align*}
    0 &= \dualsc{\mu, \pi_{\myL} v} = \scprod{h,v}_{\Ladspace}
        = \scprod{h,v}_{\Lp{2}(\Omega)^{\irow}} + \scprod{\diffopad h, \diffopad v}_{\Lp{2}(\Omega)^{\icol}}
    \\
      &= \dualprod*{h,v}_{\mathcal{D}\dual(\Omega)^{\irow},\mathcal{D}(\Omega)^{\irow}}
        +  \dualprod*{\diffopad h,{\diffopad v}}_{\mathcal{D}\dual(\Omega)^{\icol},\mathcal{D}(\Omega)^{\icol}}
    \\
      &= \dualprod*{(\idop - \diffop \diffopad) h, v}_{\mathcal{D}\dual(\Omega)^{\irow},\mathcal{D}(\Omega)^{\irow}}
      .
  \end{align*}
  This means $ \diffop \diffopad h = h$ in the sense of distributions. However, $h \in \Ladspace$ implies $h \in \Lp{2}(\Omega)$, which in turn gives $\diffop \diffopad h \in \Lp{2}(\Omega)^{\irow}$, and $\diffopad h \in \Lp{2}(\Omega)^{\icol}$. By \eqref{eq:green-id} for $g\in\Ladspace$ we have
  \begin{align*}
    \dualsc{\mu,\pi_{\myL}g} &= \scprod{h,g}_{\Ladspace}
				= \scprod{h,g}_{\Lp{2}(\Omega)^{\irow}} + \scprod{\diffopad h, \diffopad g}_{\Lp{2}(\Omega)^{\icol}}
    \\
                              &= \scprod{(\idop - \diffop \diffopad) h,g}_{\Lp{2}(\Omega)^{\irow}} + \dualsc{\Lnu \diffopad h, \pi_{\myL} g}
    \\
                              &= \dualsc[\Big]{\Lnu \underbrace{(\diffopad h)}_{=:f}, \pi_{\myL} g}
                              .
  \end{align*}
  We define $f := \diffopad h \in \Lp{2}(\Omega)^{\icol}$, which gives us $\diffop f = \diffop \diffopad h = h \in \Lp{2}(\Omega)$ and consequently $f\in \Lspace$. Hence, $\Lnu f = \mu$ completes the proof.
\end{Proof}

By \autoref{re:Lnu-equals-LnuGamma} also $\Lnu[\Gamma_{1}]: \Lspace \to \PLspace[\Gamma_{1}]\dual$ is linear, bounded and onto.

\begin{Proposition}
  $(\PLspace[\Gamma_1], \Lp{2}_{\pi}(\Gamma_1), \PLspace[\Gamma_1]\dual)$ is a quasi Gelfand triple.
\end{Proposition}

\begin{Proof}
  Let $\tilde{D}_{+} := \ran \pi_{\myL}\big\vert_{H_{\Gamma_0}^{1}(\Omega)^{\irow}}$ and let $D_{-}$ denote the corresponding set (\myref{def:Dminus}{Definition}). Then
  $\ran \indicator_{\Gamma_1} \Lnub \boundtr \subseteq D_{-}$, which is dense in $\Lp{2}_{\pi}(\Gamma_1)$ and by \myref{th:Lnu-maps-into-dualspace}{Proposition} and \myref{th:Lnu-onto}{Theorem} also dense in $\PLspace[\Gamma_1]\dual$. Hence, assertion \textup{(iv)} of \myref{prop:iota-closed-equivalences}{Proposition} is satisfied, which yields that the completion of $\tilde{D}_{+}$ and $D_{-}$ is a quasi Gelfand triple with pivot space $\Lp{2}_{\pi}(\Gamma_1)$.
\end{Proof}

\begin{Lemma}\label{le:ker-eq}
  $\ker \pi_{\myL} = \ker \Lnuad$.
\end{Lemma}

\begin{Proof}
  The following equivalences prove the statement
  \begin{align*}
    g \in \ker \pi_{\myL} &\Leftrightarrow \scprod{\pi_{\myL} g, \psi} = 0\quad \text{for all}\quad \psi \in \PLspace\dual \\
                          &\Leftrightarrow \scprod{\pi_{\myL} g, \Lnu f} = 0\quad \text{for all}\quad f \in \Lspace \\
                          &\Leftrightarrow \scprod{\Lnuad g, \pi_{\myL\hermitian} f} = 0\quad \text{for all}\quad f \in \Lspace \\
                          &\Leftrightarrow \scprod{\Lnuad g, \phi} = 0\quad \text{for all}\quad \phi \in \PLadspace \\
    &\Leftrightarrow g \in \ker \Lnuad.
    \qedhere
  \end{align*}
\end{Proof}

\section{Port Hamiltonian Systems}

In this section we will introduce port-Hamiltonian systems on multidimensional spatial domains and formulate boundary conditions which justify existence and uniqueness of solutions.
Moreover, we will parameterize all boundary conditions that provide solutions that are non-increasing in the Hamiltonian.

\begin{Definition}\label{def:PH-system}
  Let $m \in \N$ and $\myP = (P_{i})_{i=1}^{n}$, where $P_{i}$ is a Hermitian $m \times m$ matrix. Moreover, 
let $\hamiltonian : \Omega \to \K^{m\times m}$ be such that $\hamiltonian(\zeta)\hermitian = \hamiltonian(\zeta)$ and $c\idop \leq \hamiltonian(\zeta) \leq C \idop$ for a.e. $\zeta \in \Omega$ and some constants $c,C\in \R_{+}$ independent of $\zeta$.
Then we endow the space $\XH := \Lp{2}(\Omega)^{m}$ with the scalar product
\begin{align*}
  \scprod{f,g}_{\XH} := \frac{1}{2}\scprod{\hamiltonian f, g}_{\Lp{2}(\Omega)^{m}}
  = \frac{1}{2}\int_{\Omega} \scprod{\hamiltonian(\zeta) f(\zeta), g(\zeta)}_{\K^{m}} \dx[\la](\zeta)
  .
\end{align*}
Furthermore, let $P_{0} \in \K^{m\times m}$ be such that $P_{0}\hermitian = -P_{0}$. Then will call the differential equation
\begin{align}\label{eq:PH-pde}
  \begin{aligned}
    \frac{\partial}{\partial t} x(t,\zeta)	
		&= \sum_{i=1}^{n} \frac{\partial}{\partial \zeta_{i}} P_{i} \big(\hamiltonian(\zeta) x(t,\zeta)\big)
		+ P_{0}\big(\hamiltonian(\zeta)x(t,\zeta)\big),
		&&
  t\in \R_{+}, \zeta \in \Omega,\\
  x(0,\zeta) &= x_0(\zeta), &&
  \zeta \in \Omega
  \end{aligned}
\end{align}
a \emph{linear, first order port-Hamiltonian system}, where $x_0 \in \Lp{2}(\Omega)^{m}$ is the initial state. The associated \emph{Hamilitonian} $\nrg: \XH \to \R_{+}\cup \{0\}$ is defined by
\begin{align*}
  \nrg(x) := \scprod{x,x}_{\XH}
  =\frac{1}{2}\int_{\Omega} \scprod{\hamiltonian(\zeta) x(\zeta), x(\zeta)}_{\K^{m}} \dx[\la](\zeta),
\end{align*}
where $\hamiltonian$ is called the \emph{Hamiltonian density}. We will refer to $\XH$ as the \emph{state space} and to its elements as \emph{state variables} or \emph{states}.
\end{Definition}

In most applications the Hamiltonian describes the energy in the state space.

By the convention of regarding a function $x: \R_{+} \times \Omega \to \K^{m}$ as $x: \R_{+} \to \Lp{2}(\Omega;\K^{m})$ by setting $x(t) = x(t,.)$, we can rewrite the PDE \eqref{eq:PH-pde} as
\begin{align*}
  \dot{x} =  \Big( \sum_{i=1}^{n} \partial_{i}P_{i} + P_{0} \Big) \hamiltonian x
  = (\Pop + P_{0})\hamiltonian x,
  \quad
  x(0) = x_0
  .
\end{align*}

We want to add the following assumptions.

\begin{Voraussetzung}
  Let $m,\irow,\icol \in \N$ such that $m = \irow + \icol$ and let $\myL = (L_{i})_{i=1}^{n}$
  such that $L_{i} \in \K^{\irow\times \icol}$.
  Then we define $\myP = (P_{i})_{i=1}^{n}$ by
  \begin{align*}
    P_{i} = \begin{bmatrix} 0 & L_{i} \\ L_{i}\hermitian & 0 \end{bmatrix}
    .
  \end{align*}
\end{Voraussetzung}

Clearly $\myP$ contains only Hermitian matrices. Moreover, we have the identities
\begin{align*}
  \Pop = \begin{bmatrix} 0 & \diffop \\ \diffopad & 0 \end{bmatrix}, \quad
  \Pnu =
    \begin{bmatrix}
      0 & \Lnu \\
      \Lnuad & 0
    \end{bmatrix}
               , \quad
               \pi_{\myP} = \begin{bmatrix}
                 \pi_{\myL} & 0 \\
                 0 & \pi_{\myL\hermitian}
               \end{bmatrix}
                     .
\end{align*}
Corresponding to
those splittings
we want to define $(\hamiltonian x)_{1}$ and $(\hamiltonian x)_{2}$, such that
\[
  \Pop\hamiltonian x = \left[\begin{matrix} \diffop (\hamiltonian x)_{2} \\ \diffopad (\hamiltonian x)_{1}\end{matrix}\right],
  \quad
  \left[\begin{matrix}0 & \Lnu \end{matrix}\right] \hamiltonian x= \Lnu(\hamiltonian x)_{2},
  \quad
  \left[\begin{matrix}\piL & 0\end{matrix}\right] \hamiltonian x = \piL (\hamiltonian x)_{1}.
\]

\begin{Theorem}
  The operator
  \begin{align*}
    A_{0} := -(\Pop + P_{0}) \hamiltonian, \quad \dom A_{0} := \hamiltonian^{-1}( \ker \Pnu)
  \end{align*}
  is closed, skew-symmetric, and densely defined on $\XH$. Its adjoint is
  \begin{align*}
    A_{0}\adjun = (\Pop + P_{0}) \hamiltonian, \quad \dom A_{0}\adjun = \hamiltonian^{-1} (\Pspace)
    .
  \end{align*}
  Let $B_{1} = \begin{bmatrix} \pi_{\myL} & 0 \end{bmatrix}\hamiltonian$, $B_{2} = \begin{bmatrix} 0 & \Lnu \end{bmatrix}\hamiltonian$ and
  $\Psi$ the duality map of $(\PLspace, \Lp{2}(\partial \Omega), \PLspace\dual)$.
  Then $(\PLspace, B_{1}, \Psi B_{2})$ is a boundary triple for $A_{0}\adjun$.
\end{Theorem}

\begin{Proof}
  Instead of consindering $A_{0}\adjun$ as the adjoint of $A_{0}$, we just take it as a symbol. We will justify that it is in fact the adjoint of $A_{0}$ later in the proof.
  
  By \myref{le:closed-operator}{Lemma} $\Pop$ is a closed operator on $\Pspace$. Since $\hamiltonian$ is continuous, it is easy to see that $A_{0}\adjun$ is closed with domain $\hamiltonian^{-1} (\Pspace)$. Let $B\adjunX{H}$ denote the adjoint of $B$ with respect to $\scprod{.,.}_{H}$ for any Hilbert space $H$. According to \myref{re:skew-symmetric}{Remark} it is easy to see that the adjoint $((\Pop+P_{0})\hamiltonian)\adjunX{\XH}$ equals $(\Pop\adjunX{\Lp{2}} + P_{0}\adjunX{\Lp{2}})\hamiltonian = - (\Pop+P_{0})\hamiltonian$ with domain $\hamiltonian^{-1}(\dom \Pop\adjunX{\Lp{2}}) \subseteq \hamiltonian^{-1}(\Pspace)$. Hence, $(A_{0}\adjun)\adjun$ is skew-symmetric on $\XH$. Since $A_{0}\adjun$ is closed, we have $(A_{0}\adjun)\adjun[2] = A_{0}\adjun$.

  Now we know that $A_{0}\adjun$ is the adjoint of a skew-symmetric operator. So we can talk about boundary triples for $A_{0}\adjun$.
  First we note that
  \begin{equation*}
  \ran \left[\begin{smallmatrix} B_{1} \\ \Psi B_{2}\end{smallmatrix}\right] = \ran \piL \times \ran \Psi\Lnu = \PLspace \times \PLspace.
  \end{equation*} 
  Since $\hamiltonian$ is self-adjoint and $P_{0}$ is skew-adjoint, we have for $x \in \dom A_{0}\adjun$
  \begin{align*}
    \scprod{A_{0}\adjun x, x}_{\XH} + \scprod{x, A_{0}\adjun x}_{\XH}
    = \scprod{\Pop \hamiltonian x, \hamiltonian x} + \scprod{\hamiltonian x, \Pop \hamiltonian x}
    = 2 \Re \scprod{\Pop \hamiltonian x,\hamiltonian x}.
  \end{align*}
  The identity $\Pop = \left[\begin{smallmatrix}0 & \diffop \\ \diffopad & 0\end{smallmatrix}\right]$ and \myref{co:int-by-parts-extended}{Corollary} yield
  \begin{align*}
    2 \Re \scprod{\Pop \hamiltonian x,\hamiltonian x}
    &= 2\Re\scprod*{
      \begin{bmatrix} \diffop (\hamiltonian x)_{2} \\ \diffopad (\hamiltonian x)_{1} \end{bmatrix}, \begin{bmatrix} (\hamiltonian x)_{1} \\ (\hamiltonian x)_{2} \end{bmatrix}}
    = 2 \Re \dualprod{\Lnu (\hamiltonian x)_{2}, \piL (\hamiltonian x)_{1}}_{\PLspace\dual,\PLspace} \\
    &= 2 \Re \scprod{\Psi B_{2}x, B_{1}x}_{\PLspace}
      = \scprod{B_{1}x, \Psi B_{2}x}_{\PLspace} + \scprod{\Psi B_{2}x, B_{1}x}_{\PLspace}
      .
  \end{align*}
  The polarization identity implies that $(\PLspace, B_{1}, \Psi B_{2})$ is a boundary triple for $A_{0}\adjun$.

  By \myref{le:reconstruct-A0}{Lemma} $\dom A_{0} = \ker B_{1} \cap \ker B_{2}$, which is equal to
  \begin{align*}
    \ker B_{1} \cap \ker B_{2} = \hamiltonian^{-1}\big(
    \ker \begin{bmatrix} \pi_{\myL} & 0 \end{bmatrix}\cap
    \ker \begin{bmatrix} 0 & \Lnu \end{bmatrix}
                        \big)
    = \hamiltonian^{-1}\big( \ker \pi_{\myL} \times \ker \Lnu \big)
    .
  \end{align*}
  By \myref{le:ker-eq}{Lemma} this is equal to $\hamiltonian^{-1}(\ker \Pnu)$.
\end{Proof}

\begin{Remark}
  We can replace $(\PLspace,B_{1},\Psi B_{2})$ by $(\PLspace\dual, \Psi\adjun B_{1}, B_{2})$ in the previous theorem.
\end{Remark}

\begin{Theorem}\label{th:boundary-triple-for-diffop}
	Let $A_{0}\adjun$ be the operator from the previous theorem and 
	$\Phi$ the duality map associated to the quasi Gelfand triple $(\PLspace[\Gamma_{1}],\Lp{2}_{\pi}(\Gamma_{1}),\PLspace[\Gamma_{1}]\dual)$.
	Then we have $(\PLspace[\Gamma_{1}], \begin{bmatrix} \piL & 0 \end{bmatrix} \hamiltonian, \Phi \begin{bmatrix} 0 & \Lnu[\Gamma_{1}] \end{bmatrix} \hamiltonian)$ as a boundary triple for
	\begin{equation*}
	A := A_{0}\adjun\big\vert_{\hamiltonian^{-1}\big(\Ladspace[\Gamma_{0}] \times \Lspace\big)}
	.
	\end{equation*}
\end{Theorem}

\begin{Proof}
  Since we already have a boundary triple for $A_{0}\adjun$, we can show that $A$ is the adjoint of a skew-symmetric operator by \myref{pr:boundary-triple-facts}{Proposition} \textup{(iii)}. Hence, we have to check, whether $\left[\begin{smallmatrix} 0 & \opid \\ \opid & 0 \end{smallmatrix}\right]\mc C^{\perp} \subseteq \mc C$. For $B_{1}, B_{2}$ being the mappings from the previous theorem we have
	\begin{align*}
		\mc C &= \begin{bmatrix} B_{1} \\ \Psi B_{2} \end{bmatrix} \dom A = \PLspace[\Gamma_{1}] \times \PLspace \\
		\left[\begin{matrix} 0 & \opid \\ \opid & 0 \end{matrix}\right]\mc C^{\perp} &= \set{0} \times \PLspace[\Gamma_{1}]^{\perp}
			\subseteq \PLspace[\Gamma_{1}] \times \PLspace = \mc C
			.
	\end{align*}
	For $x,y \in \dom A$ we have
	\begin{align*}
		\scprod{B_{1} x,\Psi B_{2} y}_{\PLspace} &= \dualprod[\big]{\piL (\hamiltonian x)_{1}, \Lnu (\hamiltonian y)_{2}}_{\PLspace,\PLspace\dual}
		= \dualprod[\big]{\piL (\hamiltonian x)_{1}, \Lnu[\Gamma_{1}] (\hamiltonian y)_{2}}_{\PLspace[\Gamma_{1}], \PLspace[\Gamma_{1}]\dual}
		\\
		&= \scprod*{ \begin{bmatrix} \piL & 0 \end{bmatrix} \hamiltonian x, \Phi \begin{bmatrix} 0 & \Lnu[\Gamma_{1}]\end{bmatrix} \hamiltonian y}_{\PLspace[\Gamma_{1}]}
		,
	\end{align*}
	which yields \autoref{def:boundary-triple-equation} in \autoref{def:boundary-triple}.
	By $\ran \left[\begin{smallmatrix}  \piL & 0 \\ 0 & \Phi\Lnu[\Gamma_{1}] \end{smallmatrix}\right]\Big\vert_{\Ladspace[\Gamma_{0}] \times \Lspace} = \PLspace[\Gamma_{1}] \times \PLspace[\Gamma_{1}]$, the remaining \autoref{def:boundary-triple-surjective} is fulfilled.
\end{Proof}

The next theorem is \cite[Theorem 2.5]{kurula-zwart-wave}.

\begin{Theorem}
  Let $A_{0}$ be a skew-symmetric operator on a Hilbert space $X$ and $(\mathcal{B}, B_{1}, B_{2})$ be a boundary triple for $A_{0}\adjun$. Furthermore let $\mathcal{K}$ be a Hilbert space, $W_{B} = \begin{bmatrix} W_{1} & W_{2}\end{bmatrix}$, where $W_{1},W_{2} \in \Lb{\mathcal{B}, \mathcal{K}}$, and $A := A_{0}\adjun \big|_{\dom A}$, where $\dom A = \ker W_{B} \begin{bmatrix} B_{1} \\ B_{2} \end{bmatrix}$. If $\ran W_{1} - W_{2} \subseteq \ran W_{1} + W_{2}$ then the following assertions are equivalent.
  \begin{enumerate}[label = \textrm{\textup{(\roman*)}}]
  \item The operator $A$ generates a contractions semigroup on $X$.
  \item The operator $A$ is dissipative.
  \item The operator $W_{1}+ W_{2}$ is injective and the following operator inequality holds
    \begin{align*}
      W_{1}W_{2}\adjun + W_{2}W_{1}\adjun \geq 0.
    \end{align*}
  \end{enumerate}
\end{Theorem}

We will reformulate this theorem to fit our situation.

\begin{Corollary}\label{co:semigroup-generator}
  Let $\mathcal{K}$ be some Hilbert space and $W  = \begin{bmatrix}W_{1} & W_{2}\end{bmatrix} : \PLspace[\Gamma_{1}] \times \PLspace[\Gamma_{1}] \to \mathcal{K}$ a bounded linear mapping such that $\ran W_{1} - W_{2} \subseteq \ran W_{1} + W_{2}$.
  \begin{align*}
    D := \set*{  x \in \hamiltonian^{-1}(
    \Ladspace[\Gamma_{0}] \times \Lspace
    ) :
    W_{1} \begin{bmatrix} \pi_{\myL} & 0 \end{bmatrix} \hamiltonian x
    + W_{2}\Psi\begin{bmatrix} 0 & \Lnu \end{bmatrix} \hamiltonian x
                                   = 0}
                                   ,
  \end{align*}
  where $\Psi: \PLspace[\Gamma_{1}]\dual \to \PLspace[\Gamma_{1}]$ is the duality mapping corresponding to the quasi Gelfand triple.
  Then the following assertions are equivalent.
  \begin{enumerate}[label = \textrm{\textup{(\roman*)}}]
  \item $(\Pop + P_{0})\hamiltonian\big|_{D}$ generates a contractions semigroup.
  \item $(\Pop +P_{0})\hamiltonian\big|_{D}$ is dissipative.
  \item The operator $W_{1}+ W_{2}$ is injective and the following operator inequality holds
  \begin{align*}
      W_{1}W_{2}\adjun + W_{2}W_{1}\adjun \geq 0.
    \end{align*}
  \end{enumerate}
\end{Corollary}

\Autoref{co:semigroup-generator} already gives a parameterization via $W$ for all boundary conditions that make $(\Pop + P_{0})\hamiltonian$ a generator of a contractions semigroup. However, checking continuity for boundary operators which map into $\PLspace$ can be difficult. Hence, it would be appreciated to reduce the conditions on the boundary operators to conditions on better known spaces like $\Lp{2}(\partial \Omega)$.

So for the next theorem just imagine the quasi Gelfand triple to be \linebreak $(\PLspace,\Lp{2}_{\pi}(\partial\Omega),\PLspace\dual)$ to get more satisfying conditions.

The following result is a generalization of \cite[Theorem 2.6]{kurula-zwart-wave} for quasi Gelfand triple and also fixes some minor issues.

{
\renewcommand{\hs}{\mathcal{B}_{+}}%
\renewcommand{\hsdual}{\mathcal{B}_{-}}%
\renewcommand{\hsmid}{\mathcal{B}_{0}}%
\newcommand{\bvec}[1]{\begin{bmatrix}#1\end{bmatrix}}%
\newcommand{\bsmallvec}[1]{\begin{bsmallmatrix}#1\end{bsmallmatrix}}%

\begin{Theorem}\label{th:key-theorem}
  Let $(\hs,\hsmid,\hsdual)$ be a quasi Gelfand triple, $A_{0}$ be a skew-symmetric operator and $(\hs,B_{1},\Psi B_{2})$ be a boundary triple for $A_{0}\adjun$, where
  $\Psi$
  is the duality map of the Gelfand triple. For $V_{1}, V_{2} \in \Lb{\hsmid,\mathcal{K}}$ we define
  \begin{align*}
    D := \set[\bigg]{a \in \dom A_{0}\adjun : B_{1}a, B_{2}a \in \hsmid \mspace{10mu} \text{and} \mspace{10mu}  \begin{bmatrix}V_{1} & V_{2}\end{bmatrix} \begin{bmatrix}B_{1} \\ B_{2}\end{bmatrix} a = 0}
  \end{align*}
  and the operator $A := {A_{0}\adjun\big|_{D}}$. If
  \begin{enumerate}[itemsep = 3pt, topsep = 5pt, label = \textrm{\textup{(\roman*)}}]
  \item $\begin{bmatrix} V_{1}\big\vert_{\hsmid\cap\hs} & V_{2}\big\vert_{\hsmid\cap\hsdual}\end{bmatrix}$ is closed as an operator from $\hs \times \hsdual$ to $\hsmid$, 
  \item $\ker \begin{bmatrix}V_{1} & V_{2}\end{bmatrix}$ is dissipative as linear relation on $\hsmid$,
  \item $V_{1}V_{2}\adjun + V_{2}V_{1}\adjun \geq 0$ as operator on $\hsmid$,
  \end{enumerate}
  then $A$ is a generator of a contraction semigroup.
\end{Theorem}

\begin{Proof}
  It is sufficient to show that $A$ is closed, and $A$ and $A\adjun$ are dissipative.
  \medskip\par\noindent%
  %
  %
  \emph{Step 1. Showing that $A$ is closed and dissipative.} 
  \begin{align*}
    a \in D &\Leftrightarrow  \begin{bmatrix}B_{1} \\ B_{2}\end{bmatrix} a \in (\hsmid \times \hsmid) \cap \ker  \begin{bmatrix}V_{1} & V_{2}\end{bmatrix} \\
            &\Leftrightarrow  \begin{bmatrix}B_{1} \\ \Psi B_{2}\end{bmatrix} a \in \underbrace{\ker \begin{bmatrix}V_{1}\big|_{\hsmid \cap \hs} & V_{2}\Psi\adjun\big|_{\Psi(\hsmid \cap \hsdual)}\end{bmatrix}}_{=: \mathcal{C}}.
  \end{align*}
  We can write
  \begin{align*}
    \mathcal{C} = \set[\bigg]{\bvec{q \\ p} \in \hs \times \hs : q\in \hsmid, \exists \tilde{p} \in \hsmid p = \Psi \tilde{p}, V_{1} q + V_{2} \Psi\adjun p = 0}.
  \end{align*}
  For $\bsmallvec{q\\ p} \in \mathcal{C}$ we have
  \begin{align*}
    \Re \scprod{q,p}_{\hs} = \Re \scprod{q, \Psi \tilde{p}}_{\hs} =\Re \dualprod{q,\tilde{p}}_{\hs,\hsdual} = \Re \scprod{q, \tilde{p}}_{\hsmid} \leq 0,
  \end{align*}
  which implies the dissipativity of $A$ by \myref{pr:boundary-triple-facts}{Proposition}. Assumption \textup{(i)} implies that $\mc{C}$ is closed, which implies the closedness of $A$ by \myref{pr:boundary-triple-facts}{Proposition}.
  \medskip\par\noindent%
  %
  %
  \emph{Step 2. Showing that $A\adjun$ is dissipative.}
  By \autoref{pr:boundary-triple-facts} we can characterize the domain of $A\adjun$ by
  \begin{align*}
    d \in \dom A\adjun &\Leftrightarrow \bvec{B_{1} \\
    \Psi B_{2}} d \in \begin{bmatrix}0 & \idop \\ \idop & 0\end{bmatrix} \mathcal{C}^{\perp_{\hs^{2}}} \\
                       &\Leftrightarrow \bvec{\Psi B_{2} \\
    B_{1}} d \in \cl[\hs^{2}]{\ran \bvec{{V_{1}\big|_{\hsmid \cap \hs}}\adjunX{\hs} \\ {V_{2} \Psi\adjun\big|_{\Psi(\hsmid \cap \hsdual)}}\adjunX{\hs}}}
    .
  \end{align*}
  Note that if $P$ is a bounded and everywhere defined operator, and $Q$ is a linear relation, then $(PQ)\adjun = Q\adjun P\adjun$. Hence,
  \begin{align*}
  {V_{1}\big|_{\hsmid \cap \hs}}\adjunX{\hs} = (V_{1} \iota_{+})\adjun = \iota_{+}\adjun V_{1}\adjun = \Psi V_{1}\big\vert_{{V_{1}\adjun}^{-1} (\hsmid\cap\hsdual)}
  \end{align*}
  and
  \begin{align*}
  {V_{2} \Psi\adjun \big\vert_{\Psi(\hsmid\cap\hsdual)}}\adjunX{\hs} = (V_{2} \iota_{-} \Psi\adjun)\adjun = (\iota_{-} \Psi\adjun)\adjun V_{2}\adjun
  .
  \end{align*}
  From $(\Psi\iota_{-}\adjun)\adjun = \iota_{-} \Psi\adjun$ and $\iota_{-}\adjun \stackrel{\eqref{eq:iota-minus-adjun}}{=} \Psi\adjun \iota_{+}^{-1}$ follows $(\iota_{-}\Psi\adjun)\adjun = \cl{\Psi\iota_{-}\adjun} = \iota_{+}^{-1}$. Consequently,
  \begin{equation*}
  {V_{2} \Psi\adjun \big\vert_{\Psi(\hsmid\cap\hsdual)}}\adjunX{\hs} = \iota_{+}^{-1} V_{2}\adjun = V_{2}\adjun\big\vert_{{V_{2}\adjun}^{-1}(\hsmid\cap\hs)}.
  \end{equation*}
  Hence, for
  \begin{multline*}
    \bvec{x \\ y}\in
    {\ran \bvec{{V_{1}\big|_{\hsmid \cap \hs}}\adjunX{\hs} \\ {V_{2} \Psi\adjun\big|_{\Psi(\hsmid \cap \hsdual)}}\adjunX{\hs}}}
    \\
    =
    \set*{\bvec{\Psi V_{1}\adjun \\ V_{2}\adjun} k : k\in {V_{1}\adjun}^{-1}(\hsmid\cap\hsdual) \cap {V_{2}\adjun}^{-1}(\hsmid\cap\hs)}.
  \end{multline*}
  we have
  \begin{align*}
    \Re \scprod{x,y}_{\hs} &= \Re\scprod{\Psi V_{1}\adjun k, V_{2}\adjun k}_{\hs} = \Re\dualprod{V_{1}\adjun k, V_{2}\adjun k}_{\hsdual,\hs} = \Re \scprod{V_{1}\adjun k, V_{2}\adjun k}_{\hsmid} \\
                           &= \Re\scprod{V_{2}V_{1}\adjun k,k}_{\mathcal{K}} \geq 0
                           .
  \end{align*}
  Therefore, $-A\adjun$ is accretive and $A\adjun$ is dissipative.
\end{Proof}

\begin{Remark}
  In the \hyperref[th:key-theorem]{previous theorem}, it is possible to replace the condition of $\big[\begin{smallmatrix} V_{2}\vert_{\hsmid \cap \hs} & V_{2}\vert_{\hsmid \cap \hsdual}\end{smallmatrix}\big]$ being closed by
  \begin{equation*}
    \ker\cl{
    \left[
      \begin{matrix}
        V_{2}\big\vert_{\hsmid \cap \hs} & V_{2}\big\vert_{\hsmid \cap \hsdual}
      \end{matrix}
    \right]
    }
    =
    \cl{\ker
    \left[
      \begin{matrix}
        V_{2}\big\vert_{\hsmid \cap \hs} & V_{2}\big\vert_{\hsmid \cap \hsdual}
      \end{matrix}
    \right]
    }
    .
  \end{equation*}
  Then instead of $A$ the operator closure $\cl{A}$ is a generator of contraction semigroup.
\end{Remark}

\begin{Example}\label{ex:id-and-pos-operator}
  Let $M \in \Lb{\hsmid}$ be strictly positive. Then $V_{1} := \idop$, $V_{2} := M$ fulfill all conditions of the \hyperref[th:key-theorem]{previous theorem}.
  \begin{enumerate}[label = \textrm{\textup{(\roman*)}}]
  \item Setting $S = M^{\frac{1}{2}}$ and $T = M^{-\frac{1}{2}}$ in \autoref{co:closed-operator-from-Dp-times-Dm} yields $\big[\begin{smallmatrix} \idop\big\vert_{\hsmid\cap\hs} & M\big\vert_{\hsmid\cap\hsdual}\end{smallmatrix}\big]$ being closed.
  
  \item
  For $(x,y) \in \ker \begin{bmatrix}V_{1} & V_{2}\end{bmatrix}$ we have $x = -My$. Since $M$ is positive this yields
  \begin{align*}
    \Re \scprod{x,y}_{\hsmid} = \Re \scprod{-My,y} \leq 0.
  \end{align*}
  
  \item $V_{1}V_{2}\adjun + V_{2}V_{1}\adjun = M\adjun + M = 2\Re M \geq 0$.
  \end{enumerate}
  Moreover, \autoref{co:closed-operator-from-Dp-times-Dm} also implies $\big[\begin{smallmatrix} \idop\big\vert_{\hsmid\cap\hs} & M\big\vert_{\hsmid\cap\hsdual}\end{smallmatrix}\big]$ being surjective.
  Actually, it would have been enough, if $M \in \Lb{\hsmid}$ was boundedly invertible and accretive.
  Clearly, also $V_{1} := M$, $V_{2} := \idop$ fulfill all conditions.
\end{Example}

}

\section{Port-Hamiltonian Systems as Boundary Control Systems}\label{sec:boundary-control-system}

We introduce the notion of boundary control systems, scattering passive and impedance passive in the manner of \cite{scattering-impedance-passive}. We will show that a port-Hamiltonian system can be described as such a system.

\begin{Definition}
  A \emph{colligation} $\Xi := \Big(\Big[\begin{smallmatrix}G\\ L \\ K\end{smallmatrix}\Big];\Big[\begin{smallmatrix}\mc U \\ \mc X \\ \mc Y\end{smallmatrix}\Big]\Big)$
  consists of the three Hilbert spaces $\mc U$, $\mc X$, and $\mc Y$, and the three linear maps $G$, $L$, and $K$, with the same domain $\mc Z \subseteq X$ and with values in $\mc U$, $\mc X$, and $\mc Y$, respectively.
\end{Definition}

\begin{Definition}
  A colligation $\Xi := \Big(\Big[\begin{smallmatrix}G\\ L \\ K\end{smallmatrix}\Big];\Big[\begin{smallmatrix}\mc U \\ \mc X \\ \mc Y\end{smallmatrix}\Big]\Big)$ is an \emph{(internally well-posed) boundary control system}, if
  \begin{enumerate}[itemsep=7pt, topsep=7pt, partopsep=0pt, label=\textrm{\textup{(\roman*)}}]
    
    \item the operator $\left[\begin{smallmatrix} G \\ L \\ K \end{smallmatrix}\right]$ is closed from $\mc X$ to $\left[\begin{smallmatrix} \mc U \\ \mc X \\ \mc Y \end{smallmatrix}\right]$,
    
    \item the operator $G$ is surjective, and
    
    \item the operator $A := L\big\vert_{\ker G}$ generates a contraction semigroup on $\mc X$.
  \end{enumerate}
\end{Definition}

We think of the operators in this definition as determining a system via
\begin{align}
  \begin{split}
    u(t)&= Gx(t), \\
    \dot{x}(t) &= Lx(t), \quad x(0) = x_{0}, \\
     y(t) &= Kx(t).
  \end{split}
\end{align}
We call $\mc U$ the \emph{input space}, $\mc X$ the \emph{state space}, $\mc Y$ the \emph{output space} and $\mc Z$ the \emph{solution space}.

\begin{Definition}
  Let $\Xi = \Big(\Big[\begin{smallmatrix}G\\ L \\ K\end{smallmatrix}\Big];\Big[\begin{smallmatrix}\mc U \\ \mc X \\ \mc Y\end{smallmatrix}\Big]\Big)$ be a colligation. If $\Xi$ is a boundary control system such that
  \begin{equation}\label{eq:scattering-passive}
    2\Re \scprod{L x,x}_{\mc X} + \norm{Kx}_{\mc Y}^{2} \leq \norm{Gx}_{\mc U}^{2},
  \end{equation}
  then it is \emph{scattering passive}
  and it is \emph{scattering energy preserving} if we have equality in \eqref{eq:scattering-passive}.
  
  We say $\Xi$ is \emph{impedance passive} (\emph{energy preserving}), if
  $\tilde{\Xi} := \bigg(
  \bigg[\begin{smallmatrix}\frac{1}{\sqrt{2}}(G+K)\\ L \\ \frac{1}{\sqrt{2}}(G-K)\end{smallmatrix}\bigg];
  \bigg[\begin{smallmatrix}\mc U \\ \mc X \\ \mc Y\end{smallmatrix}\bigg]\bigg)$
  is scattering passive (energy preserving).
\end{Definition}

\newcommand{\Gph}{G_{\mathrm{p}}}
\newcommand{\Lph}{L_{\mathrm{p}}}
\newcommand{\Kph}{K_{\mathrm{p}}}
Corresponding to a \hyperref[def:PH-system]{port-Hamiltonian} system we want to introduce the following operators
\begin{align*}
  \Gph := T\begin{bmatrix} \piL & 0\end{bmatrix} \hamiltonian,
  \quad
  \Lph := (\Pop + P_{0}) \hamiltonian
  \quad\text{and}\quad
  \Kph := (T\adjun)^{-1}\begin{bmatrix} 0 & \Lnu[\Gamma_{1}]\end{bmatrix}\hamiltonian
  ,
\end{align*}
where $T \in \Lb{\Lp{2}(\Gamma_{1})}$ is invertible. By \autoref{th:trans-boundary-triple} also $\Gph$ and $\Kph$ establish a boundary triple for $\Lph$.
For simplification $T$ can be imagined to be the identity mapping. We still have $\Gamma_{0}, \Gamma_{1}$ as a \hyperref[def:splitting-with-thin-boundaries]{splitting with thin boundaries} of $\partial \Omega$.

\begin{Corollary}
  The colligation $\bigg(\bigg[\begin{smallmatrix}\Gph\\ \Lph \\ \Kph\end{smallmatrix}\bigg];\bigg[\begin{smallmatrix}T\PLspace[\Gamma_{1}] \\ \XH \\ (T\PLspace[\Gamma_{1}])\dual \end{smallmatrix}\bigg]\bigg)$ with solution space
  \[
    \mc Z = \hamiltonian^{-1}\big(\Ladspace[\Gamma_{0}] \times \Lspace\big)
  \]
  is a boundary control system.
\end{Corollary}

\begin{Proof}
  Since $\Lph$ is closed with domain $\mc Z$, and $\Gph$ and $\Kph$ are continuous with the graph norm of $\Lph$, we have $(\Gph,\Lph,\Kph)$ is closed. By construction $\Gph$ is surjective. Since $\Gph$ is one operator of a boundary triple for $\Lph$, the restriction $\Lph\big\vert_{\ker\Gph}$ is skew-adjoint and therefore a generator of a contraction semigroup.
\end{Proof}

\begin{Proposition}\label{th:phs-scattering-preserving}
  Let $R \in \Lb{T\Lp{2}_{\pi}(\Gamma_{1})}$ be strictly positive.
  Then the colligation $\Xi = \bigg(\bigg[\begin{smallmatrix}\frac{1}{\sqrt{2}}(\Gph+R\Kph)\\ \Lph \\ \frac{1}{\sqrt{2}}(\Gph-R\Kph)\end{smallmatrix}\bigg];\bigg[\begin{smallmatrix} \mc U \\ \XH \\ \mc Y \end{smallmatrix}\bigg]\bigg)$ with $\mc U = \mc Y = T\Lp{2}_{\pi}(\Gamma_{1})$ endowed with $\norm{f}_{\mc U} = \norm{f}_{\mc Y} = \norm{R^{\nicefrac{-1}{2}}f}_{\Lp{2}}$ and solution space
  \[
    \mc Z = \set{x \in \hamiltonian^{-1}(\Ladspace[\Gamma_{0}] \times \Lspace) : \Gph x, \Kph x \in T\Lp{2}_{\pi}(\Gamma_{1})}.
  \]
  is a scattering energy preserving boundary control system.
\end{Proposition}

\begin{Proof}
  Let $(x_{n},[\begin{matrix}\Gph x_{n} & \Lph x_{n} & \Kph x_{n}\end{matrix}]\trans)_{n\in\N}$ be a sequence in $[\begin{matrix}\Gph&\Lph&\Kph\end{matrix}]\trans$ that converges to $(x,[\begin{matrix}f &y&g\end{matrix}]\trans)$. Since $\Lph$ with domain $\Pspace$ is a closed operator and $\Ladspace[\Gamma_{0}] \times \Lspace$ is closed in $\Pspace$, we conclude that $x \in \hamiltonian^{-1}(\Ladspace[\Gamma_{0}] \times \Lspace)$. Hence, $\Gph x_{n}$ converges in $T\PLspace[\Gamma_{1}]$ to $\Gph x$. Since $(T\PLspace[\Gamma_{1}], \Lp{2}_{\pi}(\Gamma_{1}),(T\PLspace[\Gamma_{1}])\dual)$ is a quasi Gelfand triple $\Gph x = f$. Analogously, we conclude $\Kph x = g$. Therefore, $x \in \mc Z$ and $[\begin{matrix}\Gph&\Lph&\Kph\end{matrix}]\trans$ is closed, which implies that also $\big[\begin{matrix}\frac{1}{\sqrt{2}}(\Gph + R\Kph)& \Lph& \frac{1}{\sqrt{2}}(\Gph - R\Kph)\end{matrix}\big]\trans$ is closed.
  
  By \autoref{ex:id-and-pos-operator} and \autoref{th:key-theorem} $\Lph\big\vert_{\ker \frac{1}{\sqrt{2}}(\Gph + R\Kph)}$ generates a contraction semigroup.
  
  \autoref{ex:id-and-pos-operator} also gives the surjectivity of $\frac{1}{\sqrt{2}}(\Gph + R\Kph)$.
  
  Since $(\PLspace, \Gph, \Psi\Kph)$ is a boundary triple for $\Lph$, we have
  \begin{align*}
    2 \Re \scprod{\Lph x, x}_{\XH} &= 2 \Re \dualprod{\Gph x, \Kph x}_{\PLspace, \PLspace\dual} = 2 \Re \scprod{\Gph x, \Kph x}_{\Lp{2}_{\pi}(\Gamma_{1})} \\
   & = \frac{1}{2} \big(\scprod{R^{-1}\Gph x, \Gph x}_{\Lp{2}} + 2\Re\scprod{\Gph x, \Kph x}_{\Lp{2}} + \scprod{R\Kph x, \Kph x}_{\Lp{2}}\big) \\
   & \phantom{ = } \mspace{10mu} - \frac{1}{2} \big(\scprod{R^{-1}\Gph x, \Gph x}_{\Lp{2}} - 2\Re\scprod{\Gph x, \Kph x}_{\Lp{2}} + \scprod{R\Kph x, \Kph x}_{\Lp{2}}\big) \\
   & = \norm[\big]{\tfrac{1}{\sqrt{2}}(\Gph + R\Kph)x}_{\mc U}^{2} - \norm[\big]{\tfrac{1}{\sqrt{2}}(\Gph - R\Kph)x}_{\mc Y}^{2}
   ,
  \end{align*}
  which makes
  $\Xi$
  scattering energy preserving.
\end{Proof}

\begin{Remark}\label{re:phs-scattering-passive}
  Clearly, the previous theorem holds also true for the operator triple $\big[\begin{matrix}\frac{1}{\sqrt{2}}(R\Kph + \Gph) & \Lph & \frac{1}{\sqrt{2}}(R\Kph - \Gph)\end{matrix}\big]\trans$ and for $\Gph$ and $\Kph$ being swapped.
  Moreover, replacing $\Lph$ by $\Lph + J$, where $J \in \Lb{\XH}$ is dissipative, yields a scattering passive system.
\end{Remark}

So the port-Hamiltonian system with input $u$ and output $y$ described by the following equations is well-posed.
\begin{align*}
  \sqrt{2}u(t,\zeta) &= \piL\big(\hamiltonian(\zeta) x(t,\zeta)\big)_{2} + R\Lnu\big(\hamiltonian(\zeta)x(t,\zeta)\big)_{1}, && t\in \R_{+}, \zeta \in \Gamma_{1}, \\
  \frac{\partial}{\partial t} x(t,\zeta)	
    &= \sum_{i=1}^{n} \frac{\partial}{\partial \zeta_{i}} P_{i} \big(\hamiltonian(\zeta) x(t,\zeta)\big)
    + P_{0}\big(\hamiltonian(\zeta)x(t,\zeta)\big),
    &&
    t\in \R_{+}, \zeta \in \Omega,
  \\
  \sqrt{2}y(t,\zeta) &= \piL\big(\hamiltonian(\zeta) x(t,\zeta)\big)_{2} - R\Lnu\big(\hamiltonian(\zeta)x(t,\zeta)\big)_{1}, && t\in \R_{+}, \zeta \in \Gamma_{1}, \\
  0 &= \piL\big(\hamiltonian(\zeta) x(t,\zeta)\big)_{2}, &&t\in \R_{+},\zeta \in \Gamma_{0},
\end{align*}
where $\piL$ and $\Lnu$ are used a little bit sloppy. There is always a pointwise a.e. description for these mappings, but due to compact notation we use these symbols.

\begin{Corollary}\label{th:impedance-preserving}
   The colligation $\Big(\Big[\begin{smallmatrix}\Gph \\ \Lph \\ \Kph \end{smallmatrix}\Big];\Big[\begin{smallmatrix}T\Lp{2}_{\pi}(\Gamma_{1}) \\ \XH \\ T\Lp{2}_{\pi}(\Gamma_{1})\end{smallmatrix}\Big]\Big)$ with solution space
  \[
    \mc Z = \set{x \in \hamiltonian^{-1}(\Ladspace[\Gamma_{0}] \times \Lspace) : \Gph x, \Kph x \in T\Lp{2}_{\pi}(\Gamma_{1})}
  \]
  is impedance energy preserving.
\end{Corollary}

\begin{Proof}
  This is a direct consequence of \autoref{th:phs-scattering-preserving} for $R = \opid$.
\end{Proof}

\begin{Example}[Wave equation]
  Let $\rho \in \Lp{\infty}(\Omega)$ be the mass density and $T \in \Lp{\infty}(\Omega)^{n\times n}$ be the Young modulus, such that $\frac{1}{\rho} \in \Lp{\infty}(\Omega)$, $T(\zeta)\hermitian = T(\zeta)$ and $T(\zeta) \geq \delta \opid$ for a $\delta > 0$ and almost every $\zeta \in \Omega$.
  Then the wave equation
  \[
    \frac{\partial^{2}}{\partial t^{2}} w(t,\xi) = \frac{1}{\rho(\xi)} \Div \big( T(\xi) \grad w(t,\xi)\big),
  \]
  can be formulated as a port-Hamiltonian system by choosing the state variable
  $x(t,\zeta) = \left[\begin{smallmatrix} \rho(\xi) \frac{\partial}{\partial t}w(t,\zeta) \\ \grad w(t,\zeta)\end{smallmatrix}\right]$\ignorespaces
  .
  Then the PDE looks like
  \[
    \dot{x} =
    \underbrace{
    \begin{bmatrix}
      0 & \Div \\
      \grad & 0
    \end{bmatrix}
    }_{=\Pop}
    \underbrace{
    \begin{bmatrix}
      \frac{1}{\rho} & 0 \\
      0 & T
    \end{bmatrix}
    }_{=\hamiltonian} x.
  \]
  This is shown in section 3 of \cite{kurula-zwart-wave}. This is exactly the port-Hamiltonian system we get from choosing $\myL$ as in \autoref{ex:div-grad}. From \autoref{ex:div-grad-2} and \autoref{ex:div-grad-3} we know that the boundary operators are $\boundtr$ and the extension of $\nu \cdot \boundtr$. So the system
  \begin{align*}
      \sqrt{2}u(t,\zeta) &= \nu \cdot \big(T(\zeta) \grad w(t,\zeta)\big) + \frac{\partial}{\partial t}w(t,\zeta) , && t\in \R_{+}, \zeta \in \Gamma_{1}, \\
      \frac{\partial^{2}}{\partial t^{2}} w(t,\xi) &= \frac{1}{\rho(\xi)} \Div \big( T(\xi) \grad w(t,\xi)\big),
        &&
        t\in \R_{+}, \zeta \in \Omega,
      \\
      \sqrt{2}y(t,\zeta) &= \nu \cdot \big(T(\zeta) \grad w(t,\zeta)\big) - \frac{\partial}{\partial t}w(t,\zeta), && t\in \R_{+}, \zeta \in \Gamma_{1}, \\
      0 &= \frac{\partial}{\partial t}w(t,\zeta), &&t\in \R_{+},\zeta \in \Gamma_{0},
  \end{align*}
  is well-posed.
\end{Example}

\begin{Example}[Maxwell equations]
  Let $\myL = (L_{i})_{i=1}^3$ be as in \autoref{ex:rot}. In this example we have already showed $\diffop = \rot$ and $\Lnub f = \nu \times f$. The corresponding differential operator for the port-Hamiltonian PDE is
  \begin{equation*}
    \Pop =
  \begin{bmatrix}
      0 & \diffop \\
      \diffopad & 0
    \end{bmatrix}
    =
    \begin{bmatrix}
      0 & \rot \\
      -\rot & 0
    \end{bmatrix}
    .
  \end{equation*}
  We write the state as $x = \left[\begin{smallmatrix} \vec{D} \\ \vec{B}\end{smallmatrix}\right]$, where $\vec{D}, \vec{B} \in \K^{3}$. We also want to introduce the positive function $\epsilon$, $\mu$, $g$ and $r$ such that
  \[
    \epsilon,\frac{1}{\epsilon},\mu,\frac{1}{\mu}, g \in \Lp{\infty}(\Omega)
    \quad \text{and}\quad
    r,\frac{1}{r} \in \Lp{\infty}(\Gamma_{1}).
  \]
  Furthermore, we define the Hamiltonian density by $\hamiltonian(\zeta) := \left[\begin{smallmatrix} \frac{1}{\epsilon(\zeta)} & 0 \\ 0 & \frac{1}{\mu(\zeta)} \end{smallmatrix}\right]$, where each block is a $3\times 3$ matrix. At last we define $\left[\begin{smallmatrix} \vec{E} \\ \vec{H}\end{smallmatrix}\right] := \hamiltonian \left[\begin{smallmatrix} \vec{D} \\ \vec{B}\end{smallmatrix}\right]$, so that we have the same notation as in \cite{weiss-staffans-maxwell}.

  The projection on $\cl{\ran \Lnub}$ is given by $g \mapsto (\nu\times g) \times \nu$, therefore $\piL$ is the extension of $g \mapsto (\nu\times \boundtr g) \times \nu$ to $\Ladspace$. The mapping $\pi_{\tau}$ from \cite{weiss-staffans-maxwell} can be compared with $\piL$ but is not exactly the same, since they have different domains and codomains.
  
  The corresponding boundary control system is a model for the Maxwell equations that looks like
  \begin{align*}
    \sqrt{2}u(t,\zeta) &= r(\zeta)\nu(\zeta) \times \vec{H}(t,\zeta) + (\nu(\zeta)\times  \vec{E}(t,\zeta)) \times \nu(\zeta), && t\in \R_{+}, \zeta \in \Gamma_{1}, \\
    \epsilon(\zeta)\frac{\partial}{\partial t}\vec{E}(t,\zeta) &= \rot \vec{H}(t,\zeta) - g(\zeta)\vec{E}(t,\zeta), && t\in \R_{+}, \zeta \in \Omega, \\
    \mu(\zeta)\frac{\partial}{\partial t}\vec{H}(t,\zeta) &= -\rot \vec{E}(t,\zeta), && t\in \R_{+}, \zeta \in \Omega, \\
    \sqrt{2}y(t,\zeta) &= r(\zeta)\nu(\zeta) \times \vec{H}(t,\zeta) - (\nu(\zeta)\times  \vec{E}(t,\zeta)) \times \nu(\zeta), && t\in \R_{+}, \zeta \in \Gamma_{1}, \\
    0 &= (\nu(\zeta)\times  \vec{E}(t,\zeta)) \times \nu(\zeta), &&t\in \R_{+},\zeta \in \Gamma_{0},
  \end{align*}
  and is scattering passive by \autoref{re:phs-scattering-passive}, where we set $J = \left[\begin{smallmatrix} -g & 0 \\ 0 & 0\end{smallmatrix}\right]\hamiltonian$.
\end{Example}

\begin{Example}[Mindlin plate]
  Lets regard the differential operator $\Pop$ and the skew-symmetric matrix $P_{0}$ given by
  \begingroup
  \setlength\arraycolsep{4pt}
  \begin{equation*}
   \Pop :=
    \left[
    \begin{array}{ccc|ccccc}
      0 & 0 & 0 &                     0 & 0 & 0 &                                 \partial_{1} & \partial_{2} \\
      0 & 0 & 0 &                     \partial_{1} & 0 & \partial_{2} &  0 & 0 \\
      0 & 0 & 0 &                     0 & \partial_{2} & \partial_{1} & 0 & 0 \\
      \hline
      0 & \partial_{1} & 0 &      0 & 0 & 0 & 0 & 0 \\
      0 & 0 & \partial_{2} &      0 & 0 & 0 & 0 & 0 \\
      0 & \partial_{2} & \partial_{1} &   0 & 0 & 0 & 0 & 0 \\
      \partial_{1} & 0 & 0 &                0 & 0 & 0 & 0 & 0 \\
      \partial_{2} & 0 & 0 &          0 & 0 & 0 & 0 & 0 
    \end{array}
    \right]
    \mspace{-5mu},
    P_{0} :=
    \begin{bmatrix}
      0 & 0 & 0 & 0 & 0 & 0 & 0 & 0 \\
      0 & 0 & 0 & 0 & 0 & 0 & 1 & 0 \\
      0 & 0 & 0 & 0 & 0 & 0 & 0 & 1 \\
      0 & 0 & 0 & 0 & 0 & 0 & 0 & 0 \\
      0 & 0 & 0 & 0 & 0 & 0 & 0 & 0 \\
      0 & 0 & 0 & 0 & 0 & 0 & 0 & 0 \\
      0 & -1 & 0 & 0 & 0 & 0 & 0 & 0 \\
      0 & 0 & -1 & 0 & 0 & 0 & 0 & 0 
    \end{bmatrix}
    \mspace{-5mu}.
  \end{equation*}
  \endgroup
  It is easy to derive the corresponding $\myP = (P_{i})_{i=1}^{2}$ and $\myL = (L_{i})_{i=1}^{2}$. We define a Hamiltonian density by
  %
  %
  \begin{equation*}
    \hamiltonian =
    \left[
    \begin{array}{cccccccc}
      \frac{1}{\rho h} & 0 & 0 & 0 & 0 & 0 & 0 & 0 \\
      0 & \frac{12}{\rho h^{3}} & 0 & 0 & 0 & 0 & 0 & 0 \\
      0 & 0 & \frac{12}{\rho h^{3}} & 0 & 0 & 0 & 0 & 0 \\
      0 & 0 & 0 & \multicolumn{3}{c}{\multirow{3}{*}{ \scalebox{2}{$\myvec{D}_{b}$}} } & 0 & 0 \\
      0 & 0 & 0 &  			  								& &	& 0 & 0 \\
      0 & 0 & 0 &  	 										& &	& 0 & 0 \\
      0 & 0 & 0 & 0 & 0 & 0 & \multicolumn{2}{c}{\multirow{2}{*}{ \scalebox{1.3}{$\myvec{D}_{s}$}} }  \\
      0 & 0 & 0 & 0 & 0 & 0 &  &  
    \end{array}
    \right]
    ,
  \end{equation*}
  where $\rho$, $h$ are strictly positive function, $\myvec{D}_{b}(\zeta)$ is a strictly positive $3\times 3$ matrix and $\myvec{D}_{s}(\zeta)$ is strictly positive $2\times 2 $ matrix, such that all conditions on $\hamiltonian$ in \autoref{def:PH-system} are satisfied. We write the state variable $x$ as
  \[
    \myvec{\alpha} := \left[\begin{matrix}\rho h v & \rho\frac{h^{3}}{12}w_{1} & \rho\frac{h^{3}}{12}w_{2} & \kappa_{1,1} & \kappa_{2,2} & \kappa_{1,2} & \gamma_{1,3} & \gamma_{2,3}\end{matrix}\right]\trans
    ,
  \]
  where we stick to the notation in \cite{mindlin-plate-andrea} except that we changed the coordinates $x$, $y$ and $z$ to $1$, $2$ and $3$. Furthermore, we have
  \[
    \vec{e} := \hamiltonian \myvec{\alpha} = \left[\begin{matrix} v & w_{1} & w_{2} & M_{1,1} & M_{2,2} & M_{1,2} & Q_{1} & Q_{2}\end{matrix}\right]\trans.
  \]
  We don't want to go into details about the physical meaning of these state variables. We just want to make it easier to translate the results into the notation of \cite{mindlin-plate-andrea}.
  So the port-Hamiltonian PDE
  \[
    \frac{\partial}{\partial t} x
    = (\Pop + P_{0})\hamiltonian x
    \quad\text{looks like}\quad
    \frac{\partial}{\partial t} \myvec{\alpha} = (\Pop + P_{0}) \vec{e}
    .
  \]
  The corresponding boundary operator is
  \begin{equation*}
    \Lnub f =
    \begin{bmatrix}
      0 & 0 & 0 & \nu_{1} & \nu_{2} \\
      \nu_{1} & 0 & \nu_{2} & 0 & 0 \\
      \nu_{2} & \nu_{1} & 0 & 0 & 0
    \end{bmatrix}
    \begin{bmatrix}
    f_{1} \\ f_{2} \\ f_{3} \\ f_{4} \\ f_{5}
    \end{bmatrix}
    =
    \begingroup
      \def\arraystretch{1.7}
      \begin{bmatrix}
        \nu \cdot \left[\begin{smallmatrix}f_{4} \\ f_{5}\end{smallmatrix}\right] \\
        \nu \cdot \left[\begin{smallmatrix}f_{1} \\ f_{3}\end{smallmatrix}\right] \\
        \nu \cdot \left[\begin{smallmatrix}f_{3} \\ f_{2}\end{smallmatrix}\right]
      \end{bmatrix}
    \endgroup
    .
  \end{equation*}
  Since $\norm{\nu(\zeta)} = 1$, at least $\nu_{1}(\zeta) \not= 0$ or $\nu_{2}(\zeta) \not = 0$. This can be used to show that $\ran \Lnub = \Lp{2}(\partial \Omega)^{3}$. Therefore, $\piL$ is the extension of the boundary trace operator $\boundtr$ to $\Ladspace$.
  
  Since there is no direct physical meaning to the boundary variables
  \[
    \begin{bmatrix}
      0 & \Lnu
    \end{bmatrix}
    \vec{e}
    =
    \begingroup
      \def\arraystretch{1.7}
      \begin{bmatrix}
        \nu \cdot \left[\begin{smallmatrix}Q_{1} \\ Q_{2}\end{smallmatrix}\right] \\
        \nu \cdot \left[\begin{smallmatrix}M_{1,1} \\ M_{1,2}\end{smallmatrix}\right] \\
        \nu \cdot \left[\begin{smallmatrix}M_{1,2} \\ M_{2,2}\end{smallmatrix}\right]
      \end{bmatrix}
    \endgroup
    \quad\text{and}\quad
    \begin{bmatrix}
      \piL & 0 
    \end{bmatrix}
    \vec{e}
    =
    \begin{bmatrix}
      v \\ w_{1} \\ w_{2}
    \end{bmatrix}
    ,
  \]
  we define $\eta := \left[\begin{smallmatrix} -\nu_{2} \\ \nu_{1} \end{smallmatrix}\right]$ and apply the unitary transformation
  $
     T =
     \left[
     \begin{smallmatrix}
     1 & 0 & 0 \\
     0 & \nu_{1} & \nu_{2} \\
     0 & -\nu_{2} & \nu_{1}
     \end{smallmatrix}
     \right]
  $
  to obtain
  \[
    \begin{bmatrix}
      Q_{\nu} \\ M_{\nu,\nu} \\ M_{\nu,\eta}
    \end{bmatrix}
    := T
    \begingroup
      \def\arraystretch{1.7}
      \begin{bmatrix}
        \nu \cdot \left[\begin{smallmatrix}Q_{1} \\ Q_{2}\end{smallmatrix}\right] \\
        \nu \cdot \left[\begin{smallmatrix}M_{1,1} \\ M_{1,2}\end{smallmatrix}\right] \\
        \nu \cdot \left[\begin{smallmatrix}M_{1,2} \\ M_{2,2}\end{smallmatrix}\right]
      \end{bmatrix}
    \endgroup
    \quad\text{and}\quad
    \begin{bmatrix}
      v \\ w_{\nu} \\ w_{\eta}
    \end{bmatrix}
    :=
    \underbrace{(T\adjun)^{-1}}_{=T}
    \begin{bmatrix}
      v \\ w_{1} \\ w_{2}
    \end{bmatrix}
    .
  \]
  Hence, by \autoref{th:impedance-preserving} the system 
  \begin{align*}
    u &= 
      \begin{bmatrix}
        Q_{\nu} & M_{\nu,\nu} & M_{\nu,\eta}
      \end{bmatrix}\trans
      ,
      && \text{on } \R_{+} \times \Gamma_{1},
    \\
    \frac{\partial}{\partial t} \myvec{\alpha}	
    &= (\Pop + P_{0})\vec{e},
    &&
    \text{on } \R_{+} \times \Omega,
      \\
      %
    y &=
      \begin{bmatrix}
        v & w_{\nu} & w_{\eta}
      \end{bmatrix}\trans
      ,
    && \text{on }\R_{+} \times \Gamma_{1}, \\
    0 &=
      \begin{bmatrix}
        v & w_{\nu} & w_{\eta}
      \end{bmatrix}\trans
      ,
      && \text{on }\R_{+} \times \Gamma_{0},
  \end{align*}
  for the Mindlin plate is impedance energy preserving.
\end{Example}

\appendix
\section{}

The next example shows that it is possible to have \autoref{def:boundary-triple-surjective} and \autoref{def:boundary-triple-equation}  of a ``boundary triple'' for an operator $A$ (\autoref{def:boundary-triple}) without $A$ being the adjoint of a skew-symmetric operator. Moreover, it shows that in this situation \autoref{le:reconstruct-A0} does not hold. This demonstrates the importance of $A$ being the adjoint of a skew-symmetric operator in the definition.

\begin{Example}\label{ex:b-triple-counter-ex}
Let $A = \left[\begin{smallmatrix} 0 & \frac{\diffd}{\diffd \xi} \\ \frac{\diffd}{\diffd \xi} & 0\end{smallmatrix}\right]$ be an operator on $\Lp{2}(0,1)^{2}$ with $\dom A = H^{1}(0,1)^{2}$.
By \autoref{re:skew-symmetric} the operator $A$ is the adjoint of a skew-symmetric operator.
Integration by parts yields
\begin{align*}
\scprod{A f,g} + \scprod{f,A g}
	&= \int_{0}^{1} \scprod*{\left[\begin{matrix} f'_{2} \\ f'_{1} \end{matrix}\right], \left[\begin{matrix} g_{1} \\ g_{2} \end{matrix}\right]} \dx[\xi]
		+ \int_{0}^{1} \scprod*{\left[\begin{matrix} f_1 \\ f_{2} \end{matrix}\right], \left[\begin{matrix} g'_{2} \\ g'_{1} \end{matrix}\right]} \dx[\xi] \\
	&= \int_{0}^{1} (f_2^{\prime} g_1 + f'_1 g_2 + f_1 g'_2  + f_2 g'_1) \dx[\xi] = f_2 g_1 \Big\vert_{0}^{1} + f_1 g_2 \Big\vert_{0}^{1} \\
	&= f_2(1) g_1(1) - f_2(0) g_1(0) + f_1(1) g_2(0) - f_1(0) g_2(0) \\
	&= \scprod*{\underbrace{\left[\begin{matrix} f_2(1) \\ -f_2(0) \end{matrix}\right]}_{B_2 f}, \underbrace{\left[\begin{matrix} g_1(1)  \\ g_1(0)\end{matrix}\right]}_{B_1 g}}
	+ \scprod*{\underbrace{\left[\begin{matrix} f_1(1) \\ f_1(0) \end{matrix}\right]}_{B_1 f}, \underbrace{\left[\begin{matrix} g_2(1)  \\ -g_2(0)\end{matrix}\right]}_{B_2 g}}
	.
\end{align*}
Defining $B_1 f := \left[\begin{smallmatrix} f_1(1) \\ f_1(0) \end{smallmatrix}\right]$ and $B_2 f := \left[\begin{smallmatrix} f_2(1) \\ -f_2(0) \end{smallmatrix}\right]$ yields
\begin{equation}\label{eq:example-boundary-triple}
\scprod{A f,g} + \scprod{f,A g} = \scprod{B_1 f, B_2 g} + \scprod{B_2 f, B_1 g} .
\end{equation}
The mapping $\left[\begin{smallmatrix}B_{1} \\ B_{2} \end{smallmatrix}\right]: \dom A \to \R^{4}$ is surjective (this can be seen by choosing $f_{1}$ and $f_{2}$ to be linear interpolations). So $(\R^{4},B_{1},B_{2})$ is a boundary triple for $A$.

\medskip\noindent
We define $\hat{A}$ as the restriction of $A$ on $H^1_{\set{1}=0}(0,1) \times H^1_{\set{0} = \set{1}}(0,1)$, where
\begin{align*}
H^1_{\set{1}=0}(0,1) &:= \set{f \in H^1(0,1) : f(1) = 0}
,
\quad \text{and} \\
H^1_{\set{0}=\set{1}}(0,1) &:= \set{f \in H^1(0,1) : f(0) = f(1) }
.
\end{align*}
Therefore, we can reformulate \eqref{eq:example-boundary-triple} for $f,g \in \dom \hat{A}$
\begin{equation*}
\scprod{\hat{A} f,g} + \scprod{f,\hat{A} g} = -f_{1}(0) g_{2}(0) + f_{2}(0) (-g_{1}(0))
\end{equation*}
By defining $F_{1} f := -f_{1}(0)$ and $F_{2} f := f_{2}(0)$ we again have that $\left[\begin{smallmatrix}F_{1} \\ F_{2} \end{smallmatrix}\right]: \dom \hat{A} \to \R^{2}$ is surjective. However $\hat{A}$ is not the adjoint of a skew-symmetric operator. If it were, then $(\R^{2}, F_{1},F_{2})$ would be a boundary triple for $\hat{A}$ and
\begin{equation*}
\hat{A}\adjun = -\hat{A}\big\vert_{\ker F_{1} \cap \ker F_{2}} = - A\big\vert_{H^{1}_{0}(0,1)^{2}} = A\adjun.
\end{equation*}
which is not true since $\hat{A}$ is certainly not dense in $A$. In fact, with the boundary triple for $A$ we get that the adjoint of $\hat{A}$ is $-A\big\vert_{H^{1}_{\set{0}=\set{1}}(0,1) \times H^{1}_{\set{0}=0}(0,1)}$.
\end{Example}


\begin{Lemma}\label{th:weak-to-strong-convergent}
  Let $(x_n)_{n\in\N}$ be a weak convergent sequence in a Hilbert space $H$ with limit $x$. Then there exists a subsequence $(x_{n(k)})_{k\in\N}$ such that
  \[
    \norm[\bigg]{\frac{1}{N}\sum_{k=1}^{N} x_{n(k)} - x} \to 0.
  \]
\end{Lemma}

\begin{Proof}
  We assume that $x = 0$. For the general result we just need to replace $x_n$ by $x_n - x$.
  
  We define the subsequence inductively: $n(1) = 1$ and for $k > 1$ we choose $n(k)$ such that
  \[
    \abs{\scprod{x_{n(k)}, x_{n(j)}}} \leq \frac{1}{k} \quad \text{for all} \quad j < k.
  \]
  This is possible, because $(x_n)_{n\in\N}$ converges weakly to $0$. Note that by the principle of uniform boundedness $\sup_{n\in\N}\norm{x_{n}} \leq C$.
  \begin{align*}
    \norm[\bigg]{\frac{1}{N}\sum_{k=1}^{N} x_{n(k)}}^{2} &= \frac{1}{N^2}\sum_{k=1}^{N}\sum_{j=1}^{N}\scprod{x_{n(k)},x_{n(j)}} \\
    &= \frac{1}{N^{2}} \sum_{k=1}^{N} \norm{x_{n(k)}}^{2} + \frac{1}{N^{2}}\sum_{j=1}^{N} \sum_{k=j+1}^{N} 2\Re\scprod{x_{n(k)},x_{n(j)}} \\
    &\leq \frac{1}{N} C^{2} +  \frac{2}{N^{2}} \sum_{j=1}^{N} \sum_{k=j+1}^{N} \frac{1}{k} \leq \frac{C^{2}}{N} + \frac{1}{N}\ln(N) \to 0
    .
    \qedhere
  \end{align*}
\end{Proof}

\section*{Acknowledgments}

I thank my Ph.D. supervisor Birgit Jacob for providing inspiring references and being available for discussions, which led in particular to the proof of \autoref{th:iota-p-iota-m-closed}. I also thank Michael Kaltenb\"ack, who made me aware of \autoref{th:TTadjun-self-adjoint}.

\providecommand{\bysame}{\leavevmode\hbox to3em{\hrulefill}\thinspace}
\providecommand{\MR}{\relax\ifhmode\unskip\space\fi MR }
\providecommand{\MRhref}[2]{%
  \href{http://www.ams.org/mathscinet-getitem?mr=#1}{#2}
}
\providecommand{\href}[2]{#2}

\end{document}